\documentclass[12pt]{amsart}
\usepackage{amsbsy,amssymb,amscd,amsfonts,latexsym,amstext,delarray,
amsmath,graphicx} 
\input xypic

\usepackage{color}

\newtheorem{thm}{Theorem}[section]
\newtheorem{prop}[thm]{Proposition}
\newtheorem{cor}[thm]{Corollary}
\newtheorem{lem}[thm]{Lemma}

\newtheorem{defn}[thm]{Definition}
\newtheorem{rem}[thm]{Remark}

\numberwithin{equation}{section}

\def\bT{{\mathbb T}}

\def\C{{\mathbb C}}

\def\N{{\mathbb N}}

\def\R{{\mathbb R}}
\def\Z{{\mathbb Z}}

\def\cD{{\mathcal D}}

\def\cF{{\mathcal F}}
\def\cG{{\mathcal G}}

\def\cI{{\mathcal I}}

\def\cP{{\mathcal P}}

\def\cS{{\mathcal S}}

\def\GL{{\rm GL}}

\def\fF{{\mathfrak F}}

\title{Cohn--Elkies functions from Gabor frames}
\author{Yuri I.~Manin and Matilde Marcolli}
\date{2022}
\address{Max Planck Institute for Mathematics, Bonn, D-53111, Germany}
\email{manin@mpim-bonn.mpg.de}
\address{Department of Mathematics and Department of Computing and Mathematical Sciences, 
California Institute of Technology, Pasadena, CA 91125, USA}
\email{matilde@caltech.edu}

\begin{document}
\maketitle
\begin{abstract}
We investigate the relation between two different mathematical problems:
the construction of bounds on sphere packing density using Cohn--Elkies functions
and the construction of Gabor frames for signal analysis. In particular, 
we present a general construction of Cohn--Elkies functions in arbitrary dimension,
obtained from an approximate Wexel--Raz dual for Gabor frames with
Gaussian window. 
\end{abstract}

\section{Introduction}

In this paper we compare two seemingly different mathematical problems, showing that
they share a deep connection: the construction of bounds on the density of sphere
packings in Euclidean spaces, and the construction of Gabor frames for signal analysis.

\smallskip

The best currently available construction of bounds on the density of sphere packings
is provided by the method introduced in \cite{CoEl}, based on the construction of (radial) 
functions that vanish at the points of the lattice (or periodic set) with specific decay conditions
and sign conditions on the function and its Fourier transform. We refer to such functions
as {\em Cohn--Elkies functions}. This method was especially successful in Viazovska's
explicit construction, using modular forms, of one such Cohn--Elkies function proving the optimality of the
$E_8$ lattice for the sphere packing problem in dimension $8$. This construction was then adapted  in \cite{CKMRV}
to prove the optimality of the Leech lattice in dimension $24$. Despite these remarkable
achievements, in general explicit geometric constructions of Cohn--Elkies functions remain elusive,
through a numerical approximation algorithm using linear programming is described in \cite{CoEl}.

\smallskip

On the other hand, Gabor frames provide systems of filters for signal analysis that
have good encoding and decoding properties, though they do not consist of
orthogonal bases, \cite{Groch}.  A Gabor system is typically constructed by applying translation
and modulation operators parameterized by the points of a lattice (or more general
sets including periodic sets) to a window function with nice properties (for instance a
Gaussian). The main question then is whether a Gabor system constructed in this
way satisfies the frame condition (hence has good properties for signal analysis). 
This property depends crucially on the lattice. A good way of analyzing Gabor frames
and properties equivalent to the frame condition is through Wexel--Raz duality, \cite{GroKop}. 
This leads to a characterization of the frame condition for a Gabor system in $\R^n$ in terms 
of an entire function in $\C^n$ that vanishes at points of the lattice and is related to
the Bargmann transform of the Wexel--Raz dual window function. 

\smallskip

These two problems share the underlying question of the construction of a function vanishing at points of a lattice,
with assigned properties in terms of the closely related Fourier and Bargmann (or
short-time Fourier) transform. In the case of Cohn--Elkies functions one usually assumes
that the function is radial, hence vanishing on spheres containing lattice points, while
in the Gabor frame problem one typically deals with functions vanishing on hyperplanes
containing lattice points, in the sense of the general construction of \cite{Papush}.
In fact, as remarked in \cite{CoEl}, the radial hypothesis in the Cohn--Elkies case is not necessary, 
and we will consider more general such functions. 

\smallskip

There is another important direct relation between these two questions. A special class of
Gabor frames, called Grassmannian frames \cite{StrohHeath}, have the property that they
minimize (over lattices) the maximal correlation between the functions in the Gabor system. 
These are frames that most closely resemble the properties of orthogonal frames. 
It turns out that the optimization problem for the construction of Grassmannian frames  is
the same as the optimization problem for lattices achieving maximal sphere packing density. 

\smallskip

Given these relations between the two questions, it is natural to ask whether one can
use techniques from Gabor frame analysis to provide a different geometric approach
to the construction of Cohn--Elkies functions. In this paper we show that this is indeed
the case and that Wexel--Raz duality for Gabor frames provides a new approach to the
construction of Cohn--Elkies functions.

\smallskip

It is important to notice here the role of lattices. In the context of the sphere packing
density problem, it is expected that lattice solutions will be only a low-dimensional
feature, with the maximal density achievable by lattices diverging from the maximal
sphere packing density in higher-dimensions. The known cases of dimensions
$1,2,3,8,24$ are the only dimensions where an explicit lattice solution is known, and
may be the only ones. Thus, focusing on the possibility of lattice solutions is clearly
very restrictive. A conjecture of Zassenhaus expects the maximal density in any
dimension to be attainable by periodic packings, that is, sphere packings with
sphere centers on periodic sets (unions of translates of lattices). It is known
that periodic packings can approximate arbitrarily well the greatest packing density. 
After discussing the case of lattices, we show in the last section of this paper
how to adapt the construction to the case of periodic sets.

\smallskip

The construction of Cohn--Elkies functions that we discuss in this paper uses a lattice $L\subset \R^n$ 
(whose dual $L^\vee$ is the lattice whose density one wants to probe), 
together with a choice of an auxiliary lattice $K\subset \R^n$ chosen so that $\Lambda=L\times K$ gives a Gabor frame for
a Gaussian window. It is in general difficult to obtain explicit constructions of Wexel--Raz dual windows for Gabor frames. Indeed, 
even for the case of a Gaussian window, we need to use an approximate dual. It is interesting to notice that in both the problem
of constructing Cohn--Elkies functions and the problem of constructing Wexel--Raz dual windows, cases where direct explicit
constructions are known involve the use of modular forms: in dimension $8$ and $24$ for the Cohn--Elkies problem (\cite{CKMRV},
\cite{Viaz}), and for the Wexel--Raz duality in dimension one (that is, for lattices in $\R^2$), where the canonical dual window 
is expressible explicitly in terms of lattice theta functions, \cite{Jan2}. 

\smallskip

In the rest of this introductory section we present these two problems in
more detail, and we recall the background material that we need for
our main construction, which we present in the following section.

\medskip
\subsection{Cohn--Elkies functions}

In \cite{CoEl}, Cohn and Elkies obtained a bound on the density of sphere packings  
in terms of radial functions with assigned decay and sign properties of the function and
its Fourier transform. Viazovska's explicit modular forms construction \cite{Viaz} 
of such a function  famously solved the sphere packing problem in dimension $8$, 
and a generalization of the same method also gave a solution in dimension 
$24$, \cite{CKMRV}.

\smallskip

\begin{defn}\label{CEfunction}
A Cohn--Elkies function of dimension $n\in \N$ and of size $\ell\in \R^*_+$ is 
a real-valued Schwartz function $f(x)$ with real valued Fourier transform, such that
\begin{enumerate}
\item $f(x)\leq 0$ for all $\| x \| \geq \ell$;
\item $(\fF f)(\xi)\geq 0$ for all $\xi \in \R^n$;
\item $(\fF f)(0)>0$.
\end{enumerate}
\end{defn}

\smallskip

Note that condition $(\fF f)(\xi)\geq 0$, with $(\fF f)$ not identically zero, implies $f(0)>0$. 

\smallskip

\begin{rem}\label{radialrem}{\rm
In \cite{CoEl} the Cohn--Elkies functions are assumed to be real-valued radial functions, $f(x)=f_0(\| x \|)$, for
all $x\in \R^n$, with $f_0\in L^2([0,\infty), r^{n-1}dr)$ satisfying a rapid decay condition.
In this case the Fourier transform is automatically real-valued and radial, by the description of
Fourier transform of radial functions as Hankel transform. 
Also in \cite{CoEl} a more general decay condition is assumed for the Cohn--Elkies functions,
weaker than the Schwartz condition we consider here, which suffices for the use of the Poisson
summation formula. In fact, the condition was further generalized in \cite{CoKu}. Here
we consider the more restrictive class of Schwartz functions, as in \cite{Viaz}, but one
can replace this hypothesis with decay conditions as in \cite{CoEl} or \cite{CoKu}. }
\end{rem}

\smallskip

A sphere packing $\cP_L$ based on a lattice $L\subset \R^n$ is a packing of
spheres $S^{n-1}$ centered at the lattice points, with sphere diameters equal to
the length $\ell_L$ of the shortest lattice vector. 
The {\em density} $\Delta_\cP$ of a sphere packing
$\cP$ is the fraction of volume occupied by spheres, hence in the case of a 
lattice packing it is given by the ratio
\begin{equation}\label{packdensity}
 \Delta_{\cP_L} = \frac{{\rm Vol}(B_1^n(0)}{| L |} \, \left(\frac{\ell_L^n}{2}\right)^n \, ,
\end{equation}  
where $| L |={\rm Vol}(\R^n/L)$ is the covolume of the lattice and 
$$ {\rm Vol}(B_1^n(0)) =\frac{\pi^{n/2}}{\Gamma(\frac{n}{2}+1)} $$
is the volume of the unit ball in $\R^n$. In the case of a periodic lattice, based on a periodic set
consisting of $N$ translations of a lattice $L$, the density is similarly described, with
$| L |$ replaced by $|L|/N$ in \eqref{packdensity}. The {\em center-density} is defined as $\delta_\cP=\Delta_\cP/{\rm Vol}(B_1^n(0))$.
Thus, for a sphere packing $\cP_L$ based on a lattice $L\subset \R^n$, the center-density  is given by
\begin{equation}\label{centerdensity}
 \delta_L =\left( \frac{\ell_L}{2} \right)^n \, \frac{1}{| L |} \, , 
\end{equation} 
with $\ell_L$ the shortest length of $L$. 

\smallskip

Theorem~3.2 of \cite{CoEl} shows that the existence of a Cohn--Elkies function of 
dimension $n\in \N$ and size $\ell\in \R^*_+$ gives a bound on the
center-density $\delta_\cP$
\begin{equation}\label{CEbound}
\delta_\cP \leq \left(\frac{\ell}{2}\right)^n \frac{f(0)}{(\fF f)(0)} \, , 
\end{equation}
for any arbitrary sphere packing $\cP$ in $\R^n$.

\smallskip

\begin{rem}\label{covol1}{\rm
Note that in the sphere packing problem, the lattice covolume $| L |$ 
is fixed and can be taken $| L |=1$. Here we leave $| L |$ written
explicitly to highlight the dependence of the construction on $| L |$.
The reader should assume that it has a fixed value.}
\end{rem}

\smallskip

\begin{defn}\label{specialCE}
Let $L\subset \R^n$ be a lattice with shortest length $\ell_L$.
A Cohn--Elkies function of dimension $n\in \N$ and size $\ell_L$
is {\rm special}
if in addition to the properties of Definition~\ref{CEfunction} it also satisfies
\begin{equation}\label{CEid}
 \frac{1}{| L |} = \frac{f(0)}{(\fF f)(0)}\, .
\end{equation} 
\end{defn}

\smallskip

\begin{lem}\label{CElatt}
Given a lattice $L\subset \R^n$, suppose there is an associated 
special Cohn--Elkies function of dimension $n\in \N$ and size $\ell_L$,
with $\ell_L$ the shortest length of $L$. Then the lattice $L$ realizes 
the maximal density for sphere packings in $\R^n$.
\end{lem}

\proof As in \cite{CoEl}, from the Poisson summation formula 
$$ \sum_{\lambda\in L} f(x+\lambda) =\frac{1}{| L |} \sum_{\lambda'\in L^\vee} e^{-2\pi i \langle x,\lambda'\rangle} (\fF f)(\lambda')\, , $$
with $L^\vee$ the dual lattice, one obtains that 
$$ \sum_{\lambda\in L} f(\lambda) \leq f(0) $$
since each term with $\lambda\neq 0$ in the sum is non-positive, as
$\ell_L$ is the shortest length of $L$. On the other hand 
$$ \frac{1}{| L |} \sum_{\lambda'\in L^\vee}  (\fF f)(\lambda') \geq \frac{1}{| L |} (\fF f)(0)\, , $$
as all the other terms are non-negative. Thus, we have
$$ f(0)-\frac{1}{| L |} (\fF f)(0) \geq 0 $$
which gives the estimate
$$ \frac{1}{| L |} \leq \frac{f(0)}{(\fF f)(0)}\, . $$
The lattice packing is optimal if it achieves the Cohn-Elkies bound
$$ \delta_L = \left(\frac{\ell_L}{2}\right)^n \frac{f(0)}{(\fF f)(0)} $$
determined by the Cohn--Elkies function, hence if the above inequality
is optimized.
\endproof

\smallskip

We also recall the following observation from \cite{CoEl}.

\begin{cor}\label{CEzeros}
Given a lattice $L\subset \R^n$ with shortest length $\ell_L$ and covolume $|L|$, a
special Cohn--Elkies function of dimension $n\in \N$ and size $\ell_L$ 
vanishes on all the
nonzero vectors of $L$ and its Fourier transform vanishes on all the nonzero vectors
of the dual lattice $L^\vee$.
\end{cor}

\proof Since \eqref{CEid} holds, the Poisson summation formula gives
$$ \sum_{\lambda\in L\smallsetminus \{ 0 \}} f(\lambda) = \frac{1}{| L |} 
\sum_{\lambda'\in L^\vee \smallsetminus \{ 0 \}}  (\fF f)(\lambda')\, , $$
but on the left-hand-side all the terms are non-positive while on the
right-hand-side all the terms are non-negative, hence all terms vanish.
\endproof

\smallskip

We have formulated here Lemma~\ref{CElatt} and Corollary~\ref{CEzeros}
in the lattice case. For the analogous formulation in the case of periodic sets 
see \cite{CoEl}.

\medskip
\subsection{Gabor frames}\label{GaborSec}

The construction of good frames is a fundamental question in signal analysis. Unlike
orthogonal bases, frames are overdetermined and have some amount of redundancy, 
but they also have important properties,
such as optimization of the uncertainty principle (localization in both position and
frequency variables). The frame condition ensures good encoding (via the frame operator)
and decoding properties. In particular, we focus here on Gabor
frames, obtained by acting on a window function via translation and modulation
operators. The crucial question of when a Gabor system obtained by translation
and modulation of a window function satisfies the frame condition is completely
understood in the case of Gabor frames in $L^2(\R)$ with lattices $\Lambda\subset \R^2$,
while a full characterization in higher dimensions remains a more complicated problem. 

\smallskip

\begin{defn}\label{GaborFramesDef}
\begin{enumerate}
\item Given a window function $\phi$ in $L^2(\R^n)$ and a lattice $\Lambda\subset \R^{2n}$,
the Gabor system $\cG(\phi,\Lambda)=\{ \pi_\lambda \phi \}_{\lambda\in \Lambda}$ consists
of the collection of functions
\begin{equation}\label{pilambda}
 \pi_\lambda \phi(x)=e^{2\pi i \langle \eta, x \rangle} \phi(x-\xi) \, , 
\end{equation} 
for $\lambda=(\xi,\eta)\in \Lambda$. 
\item The Gabor system $\cG(\phi,\Lambda)$ is a frame (satisfies
the frame condition) if there are constants $C,C'>0$ such that, for all
$f\in L^2(\R^d)$
\begin{equation}\label{GaborFrame}
C\, \| f \|_{L^2(\R^d)} \leq \sum_{\lambda\in \Lambda} \left| \langle f, \pi_\lambda \phi \rangle \right|^2
\leq C'\, \| f \|_{L^2(\R^d)} \, . 
\end{equation}
\item The Gabor system $\cG(\phi,\Lambda)$ is a {\rm Bessel sequence} if the upper inequality of \eqref{GaborFrame} holds,
$$ \sum_{\lambda\in \Lambda} |\langle f, \pi_\lambda \phi\rangle |^2 \leq C' \| f \|^2 $$
for all $f\in L^2(\R^n)$. 
\end{enumerate}
\end{defn}

\smallskip

The frame operator $\cS=\cS_{\phi,\Lambda}$ associated to the Gabor system $\cG(\phi,\Lambda)$ is
given by
\begin{equation}\label{frameOp}
 \cS f = \sum_{\lambda\in \Lambda}  \langle f, \pi_\lambda \phi \rangle \, \pi_\lambda \phi \, .
\end{equation} 
The Gabor system $\cG(\phi,\Lambda)$ is a frame iff $\cS_{\phi,\Lambda}$ is both bounded 
and invertible on $L^2(\R^n)$ and a Bessel sequence if it is bounded.

\medskip
\subsection{Adjoint lattice}\label{AdjLatSec}

The adjoint lattice plays a crucial role in the Wexel--Raz duality for Gabor frames and in
the equivalent characterization of the frame condition in terms of sampling and interpolation
of entire functions.

\begin{defn}\label{AdjLattice}
Given a lattice $\Lambda\subset \R^{2n}$, the adjoint lattice $\Lambda^o$ is given by
\begin{equation}\label{CommAdjLatt}
\Lambda^o=\{ \lambda'\in \R^{2n}\,|\, \pi_\lambda \circ \pi_{\lambda'}=\pi_{\lambda'} \circ \pi_\lambda \, , \ \forall \lambda\in \Lambda \}\, . 
\end{equation}
with the translation-modulation operators $\pi_\lambda$ as in \eqref{pilambda}.
\end{defn}

\smallskip

We have the following equivalent description of the adjoint lattice (see Lemma~4.3.3 of \cite{GroKop}).

\begin{lem}\label{JAdjLatt}
For $\Lambda=A\Z^{2n}$ with $A\in \GL_{2n}(\R)$, the adjoint lattice
is given by 
\begin{equation}\label{adjointLatt}
\Lambda^o = J^{-1} \, (A^t)^{-1} \, \Z^{2n} 
\end{equation}
with
\begin{equation}\label{Jmatrix}
 J =\begin{pmatrix} 0 & I_n \\ -I_n & 0 \end{pmatrix}\, . 
\end{equation} 
\end{lem}

\proof
This simply follows from the fact that, for $\lambda=(\lambda_1,\lambda_2)$ and $\lambda'=(\lambda'_1,\lambda'_2)$, 
$$ \pi_{\lambda'} \circ \pi_\lambda =e^{2\pi i (\langle \lambda_1, \lambda'_2\rangle - \langle \lambda_2, \lambda'_1\rangle)}  \pi_\lambda \circ \pi_{\lambda'} $$
where the condition $1=e^{2\pi i \langle A k, J \lambda'\rangle}$ with $k\in \Z^{2n}$ holds iff 
$\langle A k, J \lambda'\rangle=\langle k, A^t J\lambda'\rangle \in \Z$ for all $k\in \Z^{2n}$,
which gives $\lambda'\in J^{-1} \, (A^t)^{-1} \, \Z^{2n}$.
\endproof

\smallskip

Note that the covolume satisfies $|\Lambda|=Vol(\R^{2n}/\Lambda)=|\det A|$, for a lattice of the form
$\Lambda =A \Z^{2n}$ for $A\in \GL_{2n}(\R)$, and for the adjoint lattice $|\Lambda^o|=|\Lambda|^{-1}$.

\smallskip

\begin{rem}\label{splitLat}{\rm
In the case of a split lattice, namely a lattice $\Lambda\subset \R^{2n}$ of the form
$\Lambda =L_1\times L_2$, with $L_1,L_2$ lattices in $\R^n$, the adjoint lattice is of the form
\begin{equation}\label{adjlatsplit}
 \Lambda^o = L_2^\vee \times L_1^\vee, 
\end{equation} 
where $L_i^\vee$ are the dual lattices of the $L_i$ in $\R^n$.}
\end{rem}

\medskip
\subsection{Wexel--Raz duality for Gabor frames}\label{WRdualSec}

The frame condition for a Gabor system $\cG(\phi,\Lambda)$ can be characterized
in terms of a duality relation, namely the existence of a dual window function $\gamma$ with
the property that the Gabor systems $\cG(\phi,\Lambda)$ and $\cG(\gamma,\Lambda)$ are
mutually orthogonal (Wexel-Raz biorthogonality relation). 

\smallskip

\begin{defn}\label{DualWindowDef}
For a Gabor system $\cG(\phi,\Lambda)$ in $L^2(\R^n)$ that is a Bessel sequence, 
a Wexel--Raz dual window $\gamma$ is a window function that satisfies the reconstruction identity
\begin{equation}\label{reconstruction}
f= \sum_{\lambda \in \Lambda} \langle f, \pi_\lambda \phi \rangle \pi_\lambda \gamma \, .
\end{equation}
\end{defn}

\smallskip

Dual windows are not unique. In particular, the {\em canonical dual window} is the
one obtained from the frame operator \eqref{frameOp} by $\gamma_{\phi,\Lambda} =\cS_{\phi,\Lambda}^{-1} \phi$.
In this case, while the frame operator \eqref{frameOp} provides the encoding
\begin{equation}\label{frameopencode}
 \cS_{\phi,\Lambda} : f \mapsto \sum_{\lambda \in \Lambda} 
\langle f, \pi_\lambda \phi \rangle \, \pi_\lambda \phi \, , 
\end{equation}
the canonical Wexel--Raz dual provides the corresponding decoding operator 
$$ \cS_{\phi,\Lambda}^{-1}: f \mapsto \sum_{\lambda \in \Lambda} \langle f, \pi_\lambda \gamma_{\phi,\Lambda} \rangle \, \pi_\lambda \gamma_{\phi,\Lambda}\, . $$

\smallskip

Dual windows can be characterized in terms of a vanishing property of their short-time Fourier transform on
the adjoint lattice.

\smallskip

\begin{defn}\label{STFTdef}
For a window function $\phi\in L^2(\R^n)$ the short-time Fourier transform of a
function $f\in L^2(\R^n)$ is given by
\begin{equation}\label{STFT}
(V_\phi f)(w) := \int_{\R^n} f(t) \bar\phi(t-u) e^{-2\pi i \langle v, t\rangle} dt =\langle f, \pi_w \phi\rangle,
\end{equation}
for $w=(u,v)\in \R^{2n}$.
\end{defn}

\smallskip

The short-time Fourier transform satisfies
$$ V_{\pi_w \phi}(\pi_w f)(z) = e^{2\pi i \langle z, J \cdot w\rangle} V_\phi f (z)\, , $$
with $J$ as in \eqref{Jmatrix}, and
$$ \langle V_\phi f, V_\gamma h \rangle_{L^2(\R^{2n}} = \langle f, h \rangle_{L^2(\R^n)} 
\overline{\langle \phi,\gamma \rangle}_{L^2(\R^n)}\, . $$
The phase factor $e^{2\pi i \langle \lambda, J \cdot z\rangle}$ satisfies, 
\begin{equation}\label{phaseJ}
e^{2\pi i \langle \lambda, J \cdot z\rangle}=1\, \forall \lambda\in \Lambda \ \Leftrightarrow \  z\in \Lambda^o\, . 
\end{equation}

\smallskip

We then have the following characterization of Wexel--Raz dual windows, see Theorem~4.4.1 of \cite{GroKop}. 

\smallskip

\begin{lem}\label{WRzerosLem}
For a Gabor system $\cG(\phi,\Lambda)$ in $L^2(\R^n)$ that is a Bessel sequence, 
a Wexel--Raz dual window $\gamma$ is a window function that satisfies 
\begin{equation}\label{biorthogonality}
\frac{1}{|\Lambda|} \langle \gamma, \pi_{\lambda'} \phi \rangle = \delta_{\lambda', 0} \, , \ \ \  \forall \lambda'\in \Lambda^o \, .
\end{equation}
\end{lem}

\proof We recall briefly the proof that \eqref{reconstruction} implies \eqref{biorthogonality}, and we refer the
reader to \cite{GroKop} for a more detailed account.
One first shows that if for two window functions $\phi$ and $\gamma$ in $L^2(\R^n)$ both Gabor systems
$\cG(\phi,\Lambda)$ and $\cG(\gamma,\Lambda)$ are Bessel sequences and 
$$ \sum_{\lambda'\in \Lambda^o} | V_\phi \gamma(\lambda') | < \infty\, , $$
then the Poisson summation formula gives
$$ \sum_{\lambda\in \Lambda} V_\phi f(z+\lambda) \overline{V_\gamma h(z+\lambda)} =\frac{1}{|\Lambda|}
\sum_{\lambda'\in \Lambda^o} V_\phi \gamma(\lambda') \overline{V_f h(\lambda')} e^{2\pi i \langle\lambda',
J z\rangle} \, , $$
for all $z\in \R^{2n}$ and for any $f,h\in L^2(\R^n)$, see Theorem~4.3.2 of \cite{GroKop}. 
For a dual window $\gamma$ one then writes 
$$ \langle f, h \rangle =\sum_{\lambda\in \Lambda} \langle \pi_z^* f, \pi_\lambda \phi \rangle \langle \pi_\lambda \gamma, \pi_z^* h \rangle = \sum_{\lambda\in \Lambda} V_\phi f(z+\lambda) \overline{V_\gamma h(z+\lambda)} $$
where the latter must be a constant function of $z\in \R^{2n}$ hence with Fourier coefficients
$$ \frac{1}{|\Lambda|} V_\phi \gamma(\lambda') \overline{V_f h(\lambda')} = \langle f, h \rangle \delta_{\lambda',0}\, , $$
see Theorem~4.4.1 of \cite{GroKop}. Thus, dual windows that satisfy \eqref{reconstruction} also  
satisfy the relation \eqref{biorthogonality}. 
\endproof 

\smallskip

The Gabor frame condition can then be equivalently formulated in terms of Wexel--Raz duality as follows
(see Theorem~4.4.1 of \cite{GroKop}).

\begin{prop}\label{WRdualFrame}
For a Gabor system $\cG(\phi,\Lambda)$ in $L^2(\R^n)$ the following properties are equivalent:
\begin{enumerate}
\item $\cG(\phi,\Lambda)$ is a frame;
\item $\cG(\phi, \Lambda^o)$ is a Bessel sequence and there is
a Wexel--Raz dual window $\gamma_{\phi,\Lambda} \in L^2(\R^n)$ (satisfying \eqref{biorthogonality})
such that $\cG(\gamma_{\phi,\Lambda},\Lambda)$ is also a Bessel sequence.
\end{enumerate}
\end{prop}

\smallskip

Thus, the problem of verifying the frame condition for Gabor systems is equivalently rephrased as the
problem of constructing Wexel--Raz dual windows satisfying the vanishing condition \eqref{biorthogonality}
on the adjoint lattice.

\medskip
\subsection{Grassmannian Gabor frames and sphere packings}\label{GrassmannSec}

In \cite{StrohHeath} a special class of frames is introduced that have the property of
minimizing correlation. Namely, frames 
$\{ \psi_\alpha \}_{\alpha \in \cI}$ such that the maximal
correlation $|\langle \psi_\alpha,\psi_\beta \rangle|$ over all $\alpha\neq \beta\in \cI$ is as
small as possible for a fixed redundancy. Such frames are called {\em Grassmannian frames}
(see Definition~\ref{GrassmFrameDef} below).

\smallskip

This question can be seen as follows; an orthonormal frame has no redundancy and
the basis elements are completely uncorrelated. Frames in general have redundancy
and for a fixed amount of redundancy this minimization problem is addressing
the question of how closely such a frame can resemble an orthonormal frame, 
in the sense of having as little correlation as possible among the basis elements.

\smallskip

\begin{rem}\label{findimGrassmann}{\rm
In a finite dimensional Hilbert space of dimension $n$, a generating set 
$\{ f_k \}_{k=1}^N$ has redundancy $\rho=N/n$. 
In this finite dimensional case, the problem of constructing Grassmannian
frames is shown in  \cite{StrohHeath} to be equivalent to the problem of
finding an arrangement of $N$ lines with largest possible angles
between them. This is in turn equivalent to constructing a spherical code with
fixed number $N$ of points and with largest possible minimal angle $\varphi$. }
\end{rem}

\smallskip

The notion of redundancy can be extended to the infinite dimensional case
in such a way that it agrees with the simple expression $\rho=N/n$ in finite
dimensions, as in Remark~\ref{findimGrassmann}.

\smallskip

\begin{defn}\label{redundancyDef}
Let $\{ z_k \}_{k\in \N}$ be a fixed choice of points in $\R^{2n}$
and let $B_k(z_k)=\{ x\in \R^n\,|\, \|z-z_k\| \leq k \}$. Given a lattice 
$\Lambda\subset \R^{2n}$, let $\Lambda_k=\Lambda\cap B_k(z_k)$.
The redundancy of $\cG(\phi,\Lambda)$ is defined as
\begin{equation}\label{redundancy}
 \rho(\phi,\Lambda):=\left( \lim_{k\to \infty} \frac{1}{\# \Lambda_k } \sum_{\lambda\in \Lambda_k} 
\langle \pi_\lambda \phi, \cS^{-1}_{\phi,\Lambda} \pi_\lambda \phi \rangle \right)^{-1} \, .
\end{equation} 
One defines $\rho^\pm(\phi,\Lambda)$ as the limsup/liminf when the limit \eqref{redundancy}
does not exist.
\end{defn}

\smallskip

It is shown in \cite{BaCaHeLa} 
that the redundancy $\rho(\phi,\Lambda)$ of a Gabor frame is equal to its
``density of label sets" $D(\phi,\Lambda)$, which is defined as
\begin{equation}\label{densitylabel}
D(\phi,\Lambda):= \lim_{k\to \infty} \frac{\# \Lambda_k}{(2k)^{2n}} \, . 
\end{equation}

\smallskip

The definition of redundancy recalled above applies to sets $\Lambda\subset \R^{2n}$
that are not necessarily lattices. In the case of lattices the notion simplifies.

\begin{rem}\label{LattRedundancy} {\rm
With all the $z_n=0$ we have the lattice covolume
$$ \lim_{k\to \infty} \frac{\# (\Lambda\cap B_k(0))}{Vol(B_{k}(0))}=\frac{1}{|\Lambda|} \, , $$
so that the redundancy is simply given by
$$ D(\phi,\Lambda)= \frac{Vol(B_1(0))}{2^{2n} \, |\Lambda|} \, , $$
where for the unit ball $B_1(0)\subset \R^{2n}$ we have $Vol(B_1(0))=\frac{\pi^n}{n!}$.
Thus, considering Gabor frames $\cG(\phi,\Lambda)$ with fixed redundancy $\rho(\phi,\Lambda)=\rho$ 
corresponds to considering lattices $\Lambda$ with fixed covolume. The Gabor frame
condition implies that the density $D(\phi,\Lambda)\geq 1$ so we can assume a fixed
covolume $|\Lambda|\leq 1$. 
}\end{rem}

\smallskip

\begin{defn}\label{GrassmFrameDef}
For a fixed window function $\phi\in L^2(\R^n)$, a Gabor frame $\cG(\phi,\Lambda)$
for a lattice $\Lambda\subset \R^{2n}$ is a Grassmannian frame if it minimizes the
maximal correlation
\begin{equation}\label{maxcorrel}
Corr(\phi,\Lambda):=\max_{\lambda\in \Lambda\smallsetminus\{ 0 \}} |\langle \phi, \pi_\lambda\phi \rangle| \, ,
\end{equation}
with the minimization taken over lattices $\Lambda$ with fixed redundancy (fixed covolume).
\end{defn}

The relation between the problem of Grassmannian Gabor frames and the problem of
lattice optimizers for the sphere packing problem can then be formulated in the following way.

\smallskip

\begin{lem}\label{packCorr}
When the window $\phi_\alpha(x)=e^{-\alpha \| x \|}$ is a Gaussian, for a lattice $\Lambda\subset \R^{2n}$
let $\Lambda_\alpha$ denote the lattice
\begin{equation}\label{Lambdaalpha}
\Lambda_\alpha:=  \{ (\lambda_1, \frac{\pi}{\alpha} \lambda_2)\in \R^{2n} \,|\, \lambda=(\lambda_1,\lambda_2)\in \Lambda \} \, .
\end{equation}
For $\lambda=(\lambda_1,\lambda_2)\in \Lambda$, we write $\lambda_\alpha$ for the corresponding
point $\lambda_\alpha=(\lambda_1,\frac{\pi}{\alpha} \lambda_2) \in \Lambda_\alpha$. 
Searching for a lattice optimizer
$$ \Lambda_{opt}={\rm arg}\, \min_\Lambda \max_{\lambda\in \Lambda\smallsetminus\{ 0 \}} 
|\langle \phi_\alpha, \pi_\lambda\phi_\alpha \rangle|
={\rm arg}\, \min_\Lambda Corr(\phi,\Lambda)  $$
is equivalent to searching for a lattice optimizer
$$ \Lambda_{opt}={\rm arg}\, \max_\Lambda \min_{\lambda\in \Lambda\smallsetminus\{ 0 \}}
 \| \lambda_\alpha \| ={\rm arg}\, \max_\Lambda \ell_{\Lambda_\alpha} \, ,  $$
maximizing the shortest length $\ell_{\Lambda_\alpha}$ for fixed covolume of $\Lambda$, 
hence for an optimizer of the sphere packing density. 
 \end{lem}
 
 \proof We can write explicitly the correlation as
 \begin{equation}\label{corrphialpha}
 \langle \phi_{\alpha}, \pi_z \phi_{\alpha} \rangle =\int_{\R^{n}} e^{-\alpha\| x \|^2} e^{-\alpha \| x-u \|^2} e^{2\pi i x\cdot v} dx 
\end{equation} 
$$ = e^{-\frac{\alpha}{2} \| u \|^2} \int_{\R^n} e^{-2\alpha \| x- u/2\|^2} e^{2\pi i x\cdot v} dx
 = e^{-\frac{\alpha}{2} \| u \|^2} e^{\pi i u\cdot v} \int_{\R^n} e^{-2\alpha \| x \|^2} e^{2\pi i x\cdot v} dx $$
 $$ = e^{\pi i u\cdot v}\,  e^{-\frac{\alpha}{2} \| u \|^2} \, \left( \frac{\pi}{2\alpha}\right)^{n/2} \, e^{-\frac{\pi^2}{2\alpha} \| v \|^2} $$
 for $z=u+iv\in \C^n$, so that we have
 \begin{equation}\label{corralpha}
  |\langle \phi_{\alpha}, \pi_z \phi_{\alpha} \rangle|=\left( \frac{\pi}{2\alpha}\right)^{n/2} \, e^{-\frac{\alpha}{2} 
  ( \| u \|^2 +\frac{\pi^2}{\alpha^2} \| v \|^2)} \, . 
 \end{equation} 
We then have 
$$ |\langle \phi_{\alpha}, \pi_\lambda \phi_{\alpha} \rangle|= \left( \frac{\pi}{2\alpha}\right)^{n/2}
e^{-\frac{\alpha}{2}\| \lambda_\alpha \|^2} \, . $$
Thus, the correlation $|\langle \phi_\alpha, \pi_\lambda\phi_\alpha \rangle|$
monotonically decreases as $\| \lambda_\alpha \|$ increases. 
Thus, the maximum is achieved on the set of shortest vectors in $\Lambda_\alpha$.
Then optimizing the lattice $\Lambda$ by making the shortest length in $\Lambda_\alpha$
as large as possible corresponds to optimizing $\Lambda$ by making the
largest correlation $|\langle \phi_\alpha, \pi_\lambda\phi_\alpha \rangle|$ over the
shortest length vectors as small as possible. 
\endproof

\medskip
\section{Cohn--Elkies functions from Wexel--Raz duality}\label{CEandWRsec}

Using an approximate construction of a Wexel--Raz dual window for 
a Gabor system with split lattice and Gaussian window, we 
obtain a general construction of Cohn--Elkies functions associated to
critical lattices in $\R^n$.

\medskip
\subsection{Approximation of Wexel--Raz dual}

We show that, given a Gaussian window function $\phi_\alpha(x)=e^{-\alpha \| x \|^2}$, such that $\cG(\phi_\alpha,\Lambda)$
is a Gabor frame, there is a (non-canonical) dual window $\gamma$ that is well approximated by
a superposition of shifted copies of $\phi_\alpha$. 

\smallskip

We first recall briefly an argument given in Theorems 1 and 2 of \cite{ChrKim}, which our statement
in Proposition~\ref{approxdual} below generalizes.

\smallskip

For a matrix $C\in \GL_n(\R)$ and a function $f\in L^2(\R^n)$ the dilation of $f$ by $C$ is defined as
\begin{equation}\label{dilation}
 (D_C f)(x):=|\det(C)|^{1/2} f(C x) \, . 
\end{equation} 
Let $\phi$ be a compactly supported real-valued function in $\R^n$,
with ${\rm supp}(\phi)\subseteq [0, N]^n$ for some $N\in \N$, 
that satisfies the partition of unity condition
\begin{equation}\label{partunity}
\sum_{k\in \Z^n} \phi(x-k) \equiv 1 \, , \ \  \forall x\in \R^n \, .
\end{equation}
By Theorems~1 and 2 of \cite{ChrKim},  if the matrices $C,B\in \GL_n(\R)$
satisfy 
\begin{equation}\label{CBnorm}
\| C^t B \| \leq \frac{1}{\sqrt{n} (2N -1) } \, ,
\end{equation}
then there is a finite subset $\cF\subset \Z^n$ and a function
\begin{equation}\label{dualtranslates}
\gamma(x) = |\det(C^t B)| \, \left( \phi(x) + 2 \sum_{k\in \cF} \phi(x+k) \right) \, ,
\end{equation}
such that the dilated functions $D_{C^{-1}}\phi$ and $D_{C^{-1}} \gamma$ generate dual
Gabor frames $$ \cG(D_{C^{-1}}\phi, \Lambda) \ \ \text{ and } \ \  \cG(D_{C^{-1}} \gamma, \Lambda). $$
This result shows that $\gamma$ is a dual window for $\cG(\phi,\Lambda)$ by showing that the biorthogonality relation
$$ \langle \gamma, \pi_\lambda \phi \rangle =\frac{1}{|\Lambda|} \delta_{\lambda,0}\, , \ \ \forall \lambda\in \Lambda^o $$
can be equivalently stated as the property that 
\begin{equation}\label{sumWRdual}
 \sum_{k\in \Z^n} \phi(x-(B^t)^{-1} n - C k) \gamma(x-Ck) =|\det(B)|\, \delta_{n,0} \, . 
\end{equation} 
Then the key properties needed to show that this relation holds are the
partition of unity relation \eqref{partunity} and the identity
\begin{equation}\label{sumWR1}
\begin{array}{c} \displaystyle{
 1=(\sum_{n\in \Gamma} \phi(x+n))^2 =\frac{1}{|\det(C^t B)|} \sum_{n\in \Gamma} \phi(x+n) \gamma (x+n)  } \\[3mm]
 \displaystyle{
= \sum_{n\in \Gamma} \phi(x+n)(\phi(x+n)+2 \sum_{\ell\in \cF} \phi(x+n+\ell))\, . }
\end{array}
\end{equation}
where $\Gamma=[0,N-1]^n \cap \Z^n$.
For the vanishing cases of \eqref{sumWRdual}, one uses the fact that for $\| B \| \leq (\sqrt{n}(2N-1))^{-1}$ 
the vanishing of \eqref{sumWRdual} for $n\neq 0$ is guaranteed by non-overlapping supports. 

\smallskip

The following statement adapts and generalizes this argument.

\smallskip

\begin{prop}\label{approxdual}
Let $\Lambda\subset \R^{2n}$ be a lattice of the form $\Lambda=L\times K$ with lattices $L,K\subset \R^n$, 
where $L=C \Z^n$ and $K=B \Z^n$ for some $C,B\in \GL_n(\R)$. Let $\Gamma_\Omega:=\Z^n \cap [-\Omega,\Omega]^n$.
Let $\phi_\alpha(x)=e^{-\alpha \| x \|^2}$ be a Gaussian window function, such that $\cG(\phi_\alpha,\Lambda)$
is a Gabor frame. For $\Omega>0$ let $\chi_\Omega$ be the characteristic function of the set $[-\Omega,\Omega]^n$ and
let $\phi_{\alpha, \Omega}(x):=\chi_\Omega(x)\, \phi_\alpha(x)$. There is a function
\begin{equation}\label{Fsymm}
 \mu: \Gamma_\Omega \to \N \ \ \text{ with } \ \  \mu_{-k}=\mu_k \, 
\end{equation} 
such that, if the matrices $B,C$ satisfy  
\begin{equation}\label{smallCB}
\| C^t B \| \leq \frac{1}{\sqrt{n} (2\Omega -1)} \, ,
\end{equation}
then the function
\begin{equation}\label{WRdualOmega}
\gamma_\Omega(x) :=|\det(C^t B)| \, (\phi_{\alpha,\Omega}(x)+ 2  \sum_{\ell \in \Gamma_\Omega}\, \mu_\ell
\, \phi_{\alpha,\Omega}(x+\ell) ) \, .
\end{equation}
is a dual frame for $\cG(\phi_{\alpha,\Omega},\Lambda)$.  Moreover, for any $\epsilon>0$ and 
an $\Omega>0$ such that $\sup |\phi_\alpha -\phi_{\alpha,\Omega}|<\epsilon$ and
$$ 2 n e^{-\frac{\pi^2}{\alpha} n (2\Omega-1)^2} <\epsilon \, , $$
if \eqref{smallCB} holds,
then the function
\begin{equation}\label{approxdualtranslates}
\gamma(x) :=|\det(C^t B)| \, (\phi_\alpha(x)+ 2  \sum_{\ell \in \Gamma_\Omega}\, \mu_\ell
\, \phi_\alpha(x+\ell) ) 
\end{equation}
satisfies the Wexel-Raz duality for $\cG(\phi_\alpha,\Lambda)$ up to an overall error of size the maximum 
between $\epsilon$ and 
$$ |\det(C^tB)| (1+2\sum \mu_\ell) \epsilon \, . $$
\end{prop}

\proof
We need to adapt the argument of \cite{ChrKim} recalled above in two ways: first to
window functions supported in a set $[-\Omega, \Omega]^n$ and then furhter
extend it from a truncated Gaussian that is compactly supported to an actual Gaussian.

We consider a window function that is compactly supported
inside $[-\Omega, \Omega]^n$. We want to showing that, in this case, the domain $\cF$ that satisfies 
\eqref{sumWR1} can be taken to be symmetric $\cF=\check\cF$. 
To this purpose, it suffices to generalize the case of Theorems~1 of \cite{ChrKim}, with $C=1$, since
the general case is then obtained as in Theorem~2 of \cite{ChrKim}. We assume that $\phi$ is a
compactly supported window function with ${\rm supp}(\phi)\subset [-\Omega,\Omega]^n$, which
satisfies the partition of unity condition \eqref{partunity}. As described above, we want to construct a
dual window $\gamma$ that satisfies the Wexel-Raz duality expressed in the form \eqref{sumWRdual}
(with $C=1$). Let $\Gamma_\Omega:= [-\Omega ,\Omega]^n \cap \Z^n$,
and let $N_\Omega:=\# \Gamma_\Omega$. Let $\Gamma_\Omega\simeq \{ n_1, \ldots, n_{N_\Omega} \}$
be a choice of an enumeration (ordering) of the set $\Gamma_\Omega$.
As in \eqref{sumWR1}, we write
$$ 1=(\sum_{n\in \Gamma_\Omega} \phi(x+n))^2 = (\phi(x+n_1)+\cdots +\phi(x+n_{N_\Omega}) \cdot (\phi(x+n_1)+\cdots +\phi(x+n_{N_\Omega}) = $$
$$ \phi(x+n_1) (\phi(x+n_1)+2\phi(x+n_2)+\cdots + 2 \phi(x+n_{N_\Omega}) $$ $$ + \phi(x+n_2)(\phi(x+n_2)+ 2\phi(x+n_3)+
\cdots 2\phi(x+n_{N_\Omega}) + \cdots $$ 
$$+ \phi(x+n_{N_\Omega}) \phi(x+n_{N_\Omega}) \, . $$
To obtain \eqref{sumWR1}, we want to first rewrite this as
$$ = \frac{1}{|\det(B)|} \sum_{j=1}^{N_\Omega} \phi(x+n_j) \cdot \gamma_{n_j}(x) \, , $$
where, as in \cite{ChrKim}, 
$$ \gamma_m(x) =|\det(B)| (\phi(x+m)+ 2 \sum_{i=1}^n \sum_{k\in E_i^m} \phi(x+k) )\, , $$
where the sets $E_i^m$ have the property that
$$ \cup_i E_i^m =\{ m'\in \Gamma_\Omega\,|\, m' > m \} $$
in the chosen ordering $\Gamma_\Omega\simeq \{ n_1, \ldots, n_{N_\Omega} \}$. 

In \cite{ChrKim} the lexicographic ordering is used on the positive quadrant $[0,\Omega]^n$, with
$$ \{ m' > m \} =\cup_i \{ m' > m \}_i := \cup_i \{ m'\,|\, m'_i > m_i \text{ and } m_j'=m_j \text{ for } i+1\leq j\leq n \}\, , $$
so that one has $E_i^m=\{ m' > m \}_i$ with
$$ E_i^m =\left\{ k\in \Z^n  \,\bigg|\, \begin{array}{ll} 0 \leq k_j \leq  \Omega  & j=1,\ldots, i-1 \\
m_i< k_i \leq \Omega  & j=i \\
k_j=m_j &  j=i+1,\ldots, n
\end{array}
\right\} $$
Moreover, one then writes
$$ \sum_{k\in E^m_i} \phi(x+k) = \sum_{k\in \cF_i} \phi(x+k+m)\, , $$
where
$$ \cF_i = \left\{ k\in \Z^n \,\bigg|\, \begin{array}{ll} | k_j | \leq  \Omega  &  j=1,\ldots, i-1 \\
1\leq k_i \leq  \Omega  & j=i \\
k_j=0 & j=i+1,\ldots, n 
\end{array}\right\}\, .  $$

Let $\epsilon=(\epsilon_1,\ldots,\epsilon_n)\in \{ \pm \}^n$ be a sequence of $n$ signs.
In $[-\Omega,  \Omega]^n$ let $Q_\epsilon$ denote the quadrant where the
$i$-th coordinate has sign $\epsilon_i$. We write $Q_+=[0, \Omega]^n$ for the
positive quadrant where all the $\epsilon_i=+$. We denote by $Y_\epsilon : Q_+\to Q_\epsilon$ the
bijection $Y_\epsilon(k)=\epsilon k := (\epsilon_i k_i)_{i=1}^n$. We have $\Gamma_\Omega=\cup_\epsilon \Gamma_\Omega^\epsilon$,
where $\Gamma_\Omega^\epsilon=\Gamma_\Omega \cap Q_\epsilon$. We identify points of $\Gamma_\Omega$ 
with pairs $(m,\epsilon)$ with $m\in \Gamma_\Omega^+=\Gamma_\Omega\cap Q_+$. We  order the set
$\{ \epsilon \}=\{ \pm \}^n$ lexicographically, with $- < +$, and we order $\Gamma_\Omega\cap Q_+$ lexicographically as in 
 \cite{ChrKim}, so that we have in $\Gamma_\Omega$
 $$ \{ (m',\epsilon') > (m,\epsilon) \} = \cup_{\epsilon',i} \{ (m',\epsilon') \,|\, \epsilon' > \epsilon \text{ or } \epsilon'=\epsilon \text{ and } m'
 > m \}   = \cup_i  \tilde E^{m,\epsilon}_i  \, , $$
 where we have
 $$  \tilde E^{m,\epsilon}_i  := \cup_{\epsilon' > \epsilon} \Gamma^{\epsilon'}_\Omega \cup E^{m,\epsilon}_i  \ \ \text{ with }
 \ \ E^{m,\epsilon}_i =Y_\epsilon(E^m_i) \, . $$
 We can then again identify the sums
$$ \sum_{k\in \tilde E^{m,\epsilon}_i} \phi(x+k) = \sum_{k\in \tilde \cF^\epsilon_i} \phi(x+k+m)\, , $$
where
$$ \tilde \cF^\epsilon_i :=  \cup_{\epsilon' > \epsilon} \Gamma^{\epsilon'}_\Omega \cup Y_\epsilon(\cF_i) \, . $$
Since there are overlaps between the sets $\tilde \cF^\epsilon_i$, the points have finite non-negative
integer multiplicities if we view the $\tilde \cF^\epsilon_i$ as subsets of $\Gamma_\Omega$, or we count
all points with multiplicity one, and consider the disjoint unions
$$ \cF^\epsilon_i :=  \sqcup_{\epsilon' > \epsilon} \Gamma^{\epsilon'}_\Omega \sqcup Y_\epsilon(\cF_i) \, . $$
Note that the set
$$ \cF:= \sqcup_{\epsilon, i} \cF^\epsilon_i $$
is invariant under the symmetry $\epsilon \leftrightarrow -\epsilon$. Thus, the corresponding set
$$ \cup_{\epsilon, i} \tilde \cF^\epsilon_i $$
consists of $\Gamma_\Omega$ with the appropriate multiplicites assigned to each of the points,
and these multiplicites are invariant with respect to the symmetry $k \mapsto -k$.
We can then write points of this set as $(k,\mu_k)$ with $k\in \Gamma_\Omega$ and $\mu_k\in \N$
the resulting multiplicity, namely the cardinality $\mu_k=\# \Pi_\Omega^{-1}(k)$ of the fiber under the projection 
$$ \Pi_\Omega: \cF= \sqcup_{\epsilon, i} \cF^\epsilon_i\twoheadrightarrow \cup_{\epsilon, i} \tilde \cF^\epsilon_i=\Gamma_\Omega \, , $$
satisfying $\mu_{-k}=\mu_k$.  Thus, we obtain the identity
$$ 1= (\sum_{n\in \Gamma_\Omega} \phi(x+n))^2 
= \frac{1}{|\det(B)|} \sum_{k\in \Z^n} \phi(x+k+n) \gamma(x+n)\, , $$
where
$$ \gamma(x)=|\det(B)| (\phi(x)+2 \sum_{m\in \Gamma_\Omega} \mu_m\,  \phi(x+m))\, .  $$

This shows that 
we can obtain in this way a Wexel-Raz dual window for $\cG(\phi_{\alpha,\Omega},\Lambda)$, 
with the truncated Gaussian window $\phi_{\alpha,\Omega}$, given by \eqref{WRdualOmega}.

We then need to further extend this result to the case of the Gaussian window $\phi_\alpha$.
While Gaussians are not compactly supported and do not satisfy the partition of
unity property, they can be well approximated by functions that satisfy both, with
arbitrarily small error. In particular, for a one-dimensional Gaussian of the form
$$ u(t)= \Delta \cdot \left( \frac{\alpha}{\pi} \right)^{1/2} e^{-\alpha t^2}\, , $$
with $\Delta>0$, the partition of unity relation \eqref{partunity} holds up to an error term
(see \cite{Bale})
\begin{equation}\label{part1error}
\sum_{k\in \Z}  u(t-k \Delta) = 1 + 2 \cos(\frac{2\pi t}{\Delta})  e^{-\frac{\pi^2}{\alpha \Delta^2}} \, .
\end{equation}
The error terms add in the case of a multidimensional Gaussian. The approximation
\eqref{approxdualtranslates} is then obtained by applying \eqref{WRdualOmega} to
a truncation $D_{C} \chi_\Omega \cdot \phi_\alpha$, where $\chi_\Omega$ is
the characteristic function of a set $[-\Omega,\Omega]^n$. 
Thus, for the Gaussian $\phi_\alpha$ we have an error term on the partition of unity
$$ |1-\sum_{k\in \Z^n} \phi_\alpha(x+ C^t B k) |\leq 2 n e^{-\frac{\pi^2}{\alpha\, \| C^t B \|^2}} \, . $$
Under the assumption that
$$ \| C^t B \| \leq \frac{1}{\sqrt{n} (2\Omega-1)} \, , $$
this error in the partition of unity relation is bounded by
$$ 2 n e^{-\frac{\pi^2}{\alpha} n (2\Omega-1)^2}\, . $$
Thus, for a given $\epsilon>0$ we can choose an $\Omega >0$ such that  both
$$ \sup |\phi_\alpha - \phi_{\alpha,\Omega} | < \epsilon\, , $$
and the error in the partition of unity relation is $2 n e^{-\frac{\pi^2}{\alpha} n (2\Omega-1)^2} <\epsilon$,
so that the window function
$$ \gamma(x):=|\det(C^t B)| (\phi_\alpha(x)+2 \sum_{m\in \Gamma_\Omega} \mu_m\,  \phi_\alpha(x+m)) $$
satisfies the Wexel-Raz duality up to an error term of size at most $\max\{ \epsilon, \det(C^tB) (1+2\sum \mu_\ell) \epsilon \}$.
\endproof

\smallskip

\begin{rem}\label{apprsize}
Given a lattice $\Lambda=L\times K$ and Gabor frames $\cG(\phi_\alpha,\Lambda)$, if $\gamma$ is
an approximate Wexel-Raz dual window constructed as above, for a chosen error size $\epsilon$ and
a corresponding cutoff of size $\Omega$, we refer to the pair $(\epsilon,\Omega)$ as the size of
the approximation. 
\end{rem}

\medskip
\subsection{Hermite constant}

\smallskip

As in \eqref{Lambdaalpha}, for a lattice $\Lambda\subset \R^{2n}$, we denote by $\Lambda_{2\sigma}$ the lattice
\begin{equation}\label{Lambdasigma}
\Lambda_{2\sigma}:=  \{ (\lambda_1, \frac{\pi}{2\sigma} \lambda_2)\in \R^{2n} \,|\, \lambda=(\lambda_1,\lambda_2)\in \Lambda \} \, .
\end{equation}
We also write $\Lambda_{2\sigma}^o:=(\Lambda^o)_{2\sigma}$ where $\Lambda^o$ is the adjoint lattice of $\Lambda$. 
(Note that this is not the same as $(\Lambda_{2\sigma})^o$, the adjoint lattice of $\Lambda_{2\sigma}$.)
We set
\begin{equation}\label{CLambdasigma}
 C_{\Lambda,\sigma} := \left(\frac{4\sigma}{\pi}\right)^{n/2}\, Corr(\phi_{2\sigma},\Lambda^o)\, , 
\end{equation} 
with $Corr(\phi_{2\sigma},\Lambda^o)$ defined as in \eqref{maxcorrel}, with a Gaussian
window $\phi_{2\sigma}(x)=e^{-2\sigma \| x \|^2}$.

\smallskip

\begin{rem}\label{criticalrem}{\rm 
The shortest length $\ell_\Lambda$ of a lattice in $\R^{2n}$ is bounded by
$$ \ell_\Lambda^2 \leq \gamma_{2n} \cdot |\Lambda|^{1/n}\, , $$
where $\gamma_{2n}$ is the Hermite constant in $\R^{2n}$. 
A lattice $\Lambda$ is {\em critical} if $\ell_\Lambda^2=\gamma_{2n} |\Lambda|^{1/n}$.
These realize the maximum lattice-packing density. }
\end{rem}

\smallskip

\begin{lem}\label{estimateHermite}
For a  lattice $\Lambda\subset \R^{2n}$ such that $\Lambda^o_{2\sigma}$ is a critical lattice we have 
\begin{equation}\label{nCLambda}
 C_{\Lambda,\sigma} \leq e^{-n \frac{\sigma\, |\Lambda^o_{2\sigma}|^{1/n}}{\pi e}} \, . 
\end{equation} 
\end{lem}

\proof As in \eqref{corrphialpha}, we have 
\begin{equation}\label{corrphi2sigma}
 \langle \phi_{2\sigma}, \pi_z \phi_{2\sigma} \rangle = e^{\pi i u\cdot v}\,  e^{-\sigma \| u \|^2} \, \left( \frac{\pi}{4\sigma}\right)^{n/2} \, e^{-\frac{\pi^2}{4\sigma} \| v \|^2}\, , 
\end{equation} 
 so that we have
 \begin{equation}\label{corr2sigma}
  |\langle \phi_{2\sigma}, \pi_z \phi_{2\sigma} \rangle|=\left( \frac{\pi}{4\sigma}\right)^{n/2} \, e^{-\sigma ( \| u \|^2 +\frac{\pi^2}{4\sigma^2} \| v \|^2)} \, . 
 \end{equation} 
 Thus, from \eqref{CLambdasigma} we obtain
$$ C_{\Lambda,\sigma} = e^{-\sigma \ell_{\Lambda^o_{2\sigma}}^2} \, . $$
 For a critical lattice we then have 
$$ C_{\Lambda,\sigma} = e^{-\sigma \gamma_{2n} |\Lambda^o_{2\sigma}|^{1/n}} \, . $$
Thus, an upper bound on $C_{\Lambda,\sigma}$ is obtained from a lower bound 
on the Hermite constant $\gamma_{2n}$.
The Minkowski--Hlawka theorem gives a lower bound for the Hermite constant of the form
$$ \gamma_{2n} \geq \left( \frac{2 \zeta(2n) }{Vol(B^{2n}_1(0))} \right)^{1/n} \, , $$
where
$$ \zeta(2n)=(-1)^{n+1} \frac{(2\pi)^{2n} B_{2n}}{2 (2n)!} \to 1 $$
for $n\to \infty$, with the Bernoulli numbers satisfying $|B_{2n}|\sim \frac{(2n)! 2}{(2\pi)^{2n}}$.
This results in a linear estimate
$$ \gamma_{2n} \geq \frac{n}{\pi e}\, . $$
 \endproof

\smallskip

We consider  lattices $\Lambda\subset \R^{2n}$ of the form $\Lambda=L\times K$ with lattices $L,K\subset \R^n$, 
where $L=C \Z^n$ and $K=B \Z^n$ for some $C,B\in \GL_n(\R)$, with $\Lambda_{2\sigma}= L \times \frac{\pi}{2\sigma} K$.
In this case we have $\Lambda_{2\sigma}^o=(\Lambda^o)_{2\sigma}=K^\vee \times \frac{\pi}{2\sigma} L^\vee$,
while $(\Lambda_{2\sigma})^o=\frac{2\sigma}{\pi} K^\vee \times L^\vee$. In this case we set
\begin{equation}\label{Corr2n}
Corr(\phi_{2\sigma},L^\vee):=\max_{\ell \in L^\vee} |\langle \phi_{2\sigma}, \pi_{i\ell} \phi_{2\sigma} \rangle|\, ,
\end{equation}
\begin{equation}\label{CLsigma}
C_{L,\sigma}:= \left(\frac{4\sigma}{\pi}\right)^{n/2}\, Corr(\phi_{2\sigma},L^\vee_{2\sigma}) \, ,
\end{equation}
with $L^\vee_{2\sigma}:=(L^\vee)_{2\sigma}=\frac{\pi}{2\sigma} L^\vee$.

\smallskip

\begin{cor}\label{LCn}
For a lattice $L\subset \R^n$ such that $L^\vee_{2\sigma}$ is a critical lattice in $\R^n$, we have
\begin{equation}\label{CLsigmaest}
 C_{L,\sigma} \leq e^{-n\, \frac{\sigma\, | L^\vee_{2\sigma} |^{1/2n}}{2\pi e}} \, .
\end{equation} 
If $L^\vee$ is a critical lattice in $\R^n$, we correspondingly have
\begin{equation}\label{CLsigmaest}
 C_{L,\sigma} \leq e^{-n\, \frac{\frac{\pi^2}{4\sigma}\, | L^\vee |^{1/2n}}{2\pi e}} \, .
\end{equation} 
\end{cor}

\proof
As in \eqref{corr2sigma} we have
\begin{equation}\label{corr2sigma2}
  |\langle \phi_{2\sigma}, \pi_{ix} \phi_{2\sigma} \rangle|=\left( \frac{\pi}{4\sigma}\right)^{n/2} \, e^{-\frac{\pi^2}{4\sigma} \| x \|^2} \, , 
\end{equation}  
with $$ \max_{\ell \in L^\vee} |\langle \phi_{2\sigma}, \pi_{i\ell} \phi_{2\sigma} \rangle|=\left( \frac{\pi}{4\sigma}\right)^{n/2} \,  
e^{-\sigma\, \ell^2_{L^\vee_{2\sigma}}}\, . $$
If we assume that $L^\vee_{2\sigma}$ is a critical lattice in $\R^n$, we have
$\ell^2_{L^\vee_{2\sigma}} = \gamma_n \cdot | L^\vee_{2\sigma} |^{1/2n}$, 
while if we assume that $L^\vee$ is a critical lattice, we have
$\ell_{L^\vee}^2=\gamma_n \cdot | L^\vee |^{1/2n}$. 
Moreover, $\ell_{L^\vee_{2\sigma}} =\frac{\pi}{2\sigma} \ell_{L^\vee}$ and $| L^\vee_{2\sigma} |=\left( \frac{\pi}{2\sigma} \right)^n | L^\vee |$
so in this case we have
$$ \ell_{L^\vee_{2\sigma}}^2 = \left( \frac{\pi}{2\sigma} \right)^2 \ell_{L^\vee}^2=\gamma_n \cdot
\left( \frac{\pi}{2\sigma} \right)^2  | L^\vee |^{1/2n} =
\left( \frac{\pi}{2\sigma} \right)^{3/2} \gamma_n \cdot | L^\vee_{2\sigma} |^{1/2n} \, . $$
\endproof

\medskip
\subsection{Preliminary estimates}

We discuss here some preliminaries for the main construction of  \S \ref{CEsec}. 

\smallskip

For a function $f\in L^2(\R^n)$ and $u\in \R^n$ we denote by $T_u f$ the translate 
\begin{equation}\label{translate}
(T_u f)(x)=f(x-u).
\end{equation}

\smallskip

Let $\phi_\alpha(x)=e^{-\alpha \| x \|^2}$, with
\begin{equation}\label{psialpha}
 \psi_\alpha(x):=\left( \frac{\pi}{\alpha} \right)^{n/2} e^{-\frac{\pi^2}{\alpha} \| x \|^2}=\fF(\phi_\alpha)\, . 
\end{equation} 

\smallskip

\begin{lem}\label{singleGaussian}
Consider a translate $g_\alpha=T_\ell \phi_\alpha$, for some $\ell\in \R^n$, and let 
$\phi_\beta(z)=e^{-\beta \| x \|^2}$ with $\psi_\beta:=\fF \phi_\beta$. Let 
\begin{equation}\label{kappa}
\kappa_{n,\Xi,\beta,\sigma} := \exp\left( -n \, \frac{\sigma \beta \Xi^{1/2n}}{2\pi e}\right) \, .
\end{equation}
For $\beta=\frac{\pi^2}{4\sigma^2}$ and $\alpha$ and $\sigma$ in the range
\begin{equation}\label{alpharangeXi}
 \alpha \leq q \cdot \pi \ \  \text{ with } \ \ 
q :=\left\{ \begin{array}{ll} \displaystyle{\frac{\sigma}{e} \Xi}  & \Xi \leq 1 \\[4mm] 
\displaystyle{\frac{\sigma}{e}} & \Xi > 1 \, ,
\end{array} \right.
\end{equation}
the following estimate holds:
\begin{equation}\label{condition}
 \int_{\R^n} g_\alpha(y) \psi_\beta ( x-y )\, dy  \geq \kappa_{n,\Xi,\beta,\sigma} \, g_\alpha(x)\, . 
\end{equation} 
\end{lem}

\proof We have
$$  \psi_\beta ( x)= \left( \frac{\pi}{\beta} \right)^{n/2} e^{-\frac{\pi^2}{\beta} \| x \|^2} $$
and for $g_\alpha(x)=e^{-\alpha \| x-\ell \|^2}$
$$ \left( \frac{\pi}{\beta} \right)^{n/2}  \int_{\R^n} g_\alpha(y) e^{-\frac{\pi^2}{\beta} \| x-y \|^2} dy =
\left( \frac{\pi}{\beta} \right)^{n/2} \int_{\R^{n}} 
e^{-(\alpha \| y \|^2+ \frac{\pi^2}{\beta} \| y-(x-\ell) \|^2)} dy $$
$$ = \left( \frac{\pi}{\beta} \right)^{n/2}  e^{-\frac{\alpha\frac{\pi^2}{\beta}}{\alpha+\frac{\pi^2}{\beta}} \| x-\ell \|^2} \int_{\R^{n}} e^{-(\alpha + \frac{\pi^2}{\beta}) \| y-\frac{\alpha}{\alpha + \frac{\pi^2}{\beta}} (x-\ell) \|^2}\, dy $$
$$ = \left( \frac{\pi}{\beta} \right)^{n/2} 
\left( \frac{\pi}{\alpha+\frac{\pi^2}{\beta}}\right)^{n/2} \, e^{-\frac{\alpha\frac{\pi^2}{\beta}}{\alpha+\frac{\pi^2}{\beta}} \| x-\ell \|^2} =
\left( \frac{\pi^2}{\beta\alpha+\pi^2}\right)^{n/2} \, \exp\left(
-\frac{\alpha\pi^2}{\beta\alpha+\pi^2} \| x- \ell \|^2 \right)  \, . $$
 In order to verify the condition
  $$ \left( \frac{\pi^2}{\beta\alpha+\pi^2}\right)^{n/2} \cdot \exp\left(
-\frac{\alpha\pi^2}{\beta\alpha+\pi^2} \| x- \ell \|^2
\right) \geq  \exp\left( -\frac{n}{2}\, \frac{\sigma\beta \, \Xi^{1/2n} }{\pi e}  \right) \cdot 
\exp\left( - \alpha \| x- \ell \|^2 \right) $$
note that 
$$ \frac{\alpha\pi^2}{\beta\alpha+\pi^2}  < \alpha $$
is always verified since $\alpha,\beta>0$ hence
$\alpha\pi^2< \alpha\pi^2+\beta\alpha^2$. 
Thus, it suffices to check when
\begin{equation}\label{checkwhen}
 \log \left( 1+ \frac{\beta\alpha}{\pi^2} \right) \leq  \frac{\sigma \beta  \, \Xi^{1/2n} }{\pi e}\, . 
\end{equation} 
We are assuming that $\beta=\frac{\pi^2}{4\sigma^2}$, so the above gives
$$ \log \left( 1+ \frac{\alpha}{4\sigma^2} \right) \leq \frac{\pi  \, \Xi^{1/2n} }{4 \sigma e} \, . $$
When $\alpha$ and $\sigma$ are in the range \eqref{alpharange} this is satisfied.
\endproof

\begin{rem}\label{XiL}{\rm
We will be interested in the case where $\Xi=| L^\vee |$ for a lattice $L\subset \R^n$ with
$$ \kappa_{n,\Xi,\beta,\sigma} := \exp\left( -n \, \frac{\sigma \beta | L^\vee |^{1/2n}}{2\pi e}\right) \, . $$
Checking \eqref{checkwhen} is then equivalent to checking when 
$$ \left( \frac{\pi^2}{\beta\alpha+\pi^2}\right)^{n/2} \geq \exp\left( -\frac{n}{2}\, \frac{\sigma \beta \, | L^\vee_{2\sigma} |^{1/2n} }{\pi e}  
\right) $$
and this is satisfied with $q=q_L$ in the range
\begin{equation}\label{alpharange}
 \alpha \leq q_L \cdot \pi \ \  \text{ with } \ \ 
q_L :=\left\{ \begin{array}{ll} \displaystyle{\frac{\sigma}{e} | L^\vee |}  & | L^\vee |\leq 1 \\[4mm] 
\displaystyle{\frac{\sigma}{e}} & | L^\vee | > 1 \, .
\end{array} \right.
\end{equation}
}\end{rem}

\medskip
\subsection{Construction of Cohn-Elkies functions}\label{CEsec}

As above, let $\phi_\alpha(x)=e^{-\alpha \| x \|^2}$ with $\psi_\alpha=\fF(\phi_{\alpha})$ as in \eqref{psialpha}.
Let $\Lambda \subset \R^{2n}$ be a lattice, with $\Lambda_{2\sigma}$ the lattice \eqref{Lambdasigma}
and $\Lambda_{2\sigma}^o=(\Lambda^o)_{2\sigma}$.
We assume that the lattice $\Lambda$ is such that the Gabor system $\cG(\phi_\alpha,\Lambda_{2\sigma})$ satisfies
the frame condition.

\smallskip

We consider in particular lattices $\Lambda=L\times K$, as above, with 
$\Lambda_{2\sigma}= L \times \frac{\pi}{2\sigma} K$ and 
$\Lambda_{2\sigma}^o=K^\vee \times \frac{\pi}{2\sigma} L^\vee$. 

\smallskip

For $\Lambda=L\times K$, let $\cF\subset L$ be a finite subset and let $\mu: \cF\to \N$ be a function that
assigns multiplicities to the points of $\cF$. Let $\cD_{\cF,\mu}$ be the function
\begin{equation}\label{DirKer}
\cD_{\cF,\mu}(x) =1+2\sum_{\ell \in \cF} \mu_\ell \, e^{2\pi i \langle \ell, x \rangle} \, .
\end{equation}

\smallskip

For example, for $n=1$ and $\cF=([-N, N]\smallsetminus \{ 0\})\cap \Z$, with all multiplicities equal to one, 
this is related to the usual Dirichlet kernel by
$$ \cD_{([-N, N]\smallsetminus \{ 0\})\cap \Z, 1}(t)=2\sum_{k=-N}^N e^{2\pi i kx}-1 =\frac{\sin((2N+1)\pi x)}{\sin(\pi x)} -1\, . $$
We write 
\begin{equation}\label{eell}
e_\ell(x):=e^{2\pi i \langle \ell, x \rangle} \, ,
\end{equation}
so that $\cD_{\cF,\mu} =1+2\sum_{\ell \in \cF} \mu_\ell\, e_\ell$.

\smallskip

For a function $f\in L^2(\R^n)$ and $u\in \R^n$ we denote by $T_u f$ the translate 
as in \eqref{translate} and we write
\begin{equation}\label{Df}
\bT_{\cF,\mu} f :=f +2\sum_{\ell \in \cF} \mu_\ell \, T_\ell f \, . 
\end{equation}

\smallskip

Consider then functions of the form
\begin{equation}\label{hCEfunct}
h_{\Lambda,\sigma}(x):=\left( \frac{4\sigma}{\pi}\right)^{n/2}  |\langle \phi_{2\sigma}, \pi_{ix} \phi_{2\sigma}\rangle | - C_{L,\sigma} 
\end{equation}
with  $C_{L,\sigma}$ as in \eqref{CLsigma}, and $\phi_\sigma(x):=e^{-\sigma \| x \|^2}$, and
\begin{equation}\label{fCEfunct}
 f_{\cF,\mu} (x) :=\langle \bT_{\cF,\mu} \gamma , \pi_{ix} \phi_\alpha \rangle \cdot h_\Lambda(x) \, .
\end{equation}
with $\gamma=\gamma_{\phi_\alpha,\Lambda_{2\sigma}}$ a Wexel--Raz dual window for $\cG(\phi_\alpha,\Lambda_{2\sigma})$.

\smallskip

\begin{thm}\label{mainthm}
Let $\Lambda=L\times K$ be a lattice in $\R^{2n}$ such that $L^\vee$ is a critical lattice in $\R^n$ and
the lattice $K\subset \R^n$ is chosen so that the Gabor system $\cG(\phi_\alpha, \Lambda_{2\sigma})$ satisfies 
the frame condition and \eqref{smallCB} is satisfied. 
Let $\gamma=\gamma_{\phi_\alpha,\Lambda_{2\sigma}}$ be the approximation to
a Wexel--Raz dual window with
approximation size  $(\epsilon,\Omega)$ (see Remark~\ref{apprsize})
constructed as in Proposition~\ref{approxdual}. Consider a datum $(\cF,\mu)$ given by the
pair $(\Gamma_\Omega,\mu)$ of Proposition~\ref{approxdual}. Then, for $\alpha$ and $\sigma$ in the
range \eqref{alpharange}, the function \eqref{fCEfunct},
$$ f_{\Gamma_\Omega,\mu} (x) :=\langle \bT_{\Gamma_\Omega,\mu} \gamma , \pi_{ix} \phi_\alpha \rangle 
\cdot h_\Lambda(x) $$
is a Cohn-Elkies function of dimension $n$ and size $\ell_{L^\vee}$, the shortest length of $L^\vee$.
\end{thm}

\smallskip

\proof We first show that $f_{\Gamma_\Omega,\mu}$ is a real valued Schwartz function that
satisfies  $f_{\Gamma_\Omega,\mu}\leq 0$ for $\| x \| \geq \ell_{L^\vee}$.
We have as in \eqref{corr2sigma2}
$$ |\langle \phi_{2\sigma}, \pi_{ix} \phi_{2\sigma}\rangle | =\left( \frac{\pi}{4\sigma}\right)^{n/2} \, e^{-\frac{\pi^2}{4\sigma} \| x \|^2}\, ,  $$
so that we have
$$ h_\Lambda(x) =\phi_{\frac{\pi^2}{4\sigma^2}}( x) - C_{L,\sigma} =\phi_{\sigma}(\frac{\pi}{2\sigma} x) - C_{L,\sigma} \, . $$
Since $C_{L,\sigma}=e^{-\sigma \ell^2_{L^\vee_{2\sigma}}}$, we have
$$ h_\Lambda(x) \leq 0 \ \ \ \text{ for } \ \ \  \| x \| \geq \ell_{L^\vee_{2\sigma}} \, \frac{2\sigma}{\pi}=\ell_{L^\vee} \, . $$

\smallskip

For $z=u+iv\in \C^n$, we have
$$ \langle \gamma, \pi_z \phi_\alpha \rangle =
\int_{\R^n} \gamma(x) e^{2\pi i x\cdot v} \phi_\alpha(x-u) \, dx =\fF( \gamma\cdot T_u \phi_\alpha) =\fF(\gamma)\star ( e_u \cdot \fF(\phi_\alpha))\, , $$
with $e_u$ as in \eqref{eell}. Thus, we have
$$ \langle \gamma, \pi_z \phi_\alpha \rangle |_{i\R^n}=\langle \gamma, \pi_{iv} \phi_\alpha \rangle =\fF( \gamma\cdot \phi_\alpha)(v) \, . $$
and we obtain 
$$  \langle \bT_{\Gamma_\Omega,\mu} \gamma, \pi_{ix} \phi_\alpha \rangle=
\langle \gamma , \pi_{ix} \phi_\alpha \rangle+2 \sum_{\ell \in \cF}  \langle T_\ell \gamma , \pi_{ix} \phi_\alpha \rangle =
\fF( \gamma)\star \fF(\phi_\alpha) +
2 \sum_{\ell \in \cF} \mu_\ell \fF(T_\ell \gamma)\star \fF(\phi_\alpha)\, . $$
We also have
$$ \fF(\bT_{\Gamma_\Omega,\mu} \gamma)= \fF( \gamma)+2\sum_{\ell \in \cF} \mu_\ell \fF(T_\ell \gamma) = \cD_{\cF,\mu} \cdot \fF(\gamma) \, . $$

\smallskip

We consider here the case where the pair $(\cF,\mu)$ is given by $(\Gamma_\Omega,\mu)$
as in \eqref{Fsymm} in Proposition~\ref{approxdual}. Since both the set $\Gamma_\Omega$ and the
multiplicity function $\mu: \Gamma_\Omega\to \N$ are invariant under the symmetry $x\mapsto -x$,
the function $\cD_{\Gamma_\Omega,\mu}$ also satisfies the symmetry 
$$ \cD_{\Gamma_\Omega,\mu}(-x)=\cD_{\Gamma_\Omega,\mu}(x) \, . $$
Since $\bar e_\ell(x)=e_\ell(-x)$, we also have $\bar \cD_{\Gamma_\Omega,\mu}(x)=\cD_{\Gamma_\Omega,\mu}(-x)$,
hence $\cD_{\Gamma_\Omega,\mu}$ is a real-valued even function. 

\smallskip

We use as $\gamma(x)$ the Wexel-Raz dual window approximation of Proposition~\ref{approxdual}, given by
\eqref{approxdualtranslates}.  We then have
$$ \fF(\gamma)= |\det(C^t B)| \, (\fF(\phi_\alpha)+ 2  \sum_{\ell \in \Gamma_\Omega}\, \mu_\ell
\, \fF(T_\ell \phi_\alpha) ) =|\det(C^t B)| \, \cD_{\Gamma_\Omega,\mu} \cdot \fF(\phi_\alpha) \, .  $$
Thus, we obtain
$$ f_{\Gamma_\Omega,\mu} =|\det(C^t B)| \,  \cD^2_{\Gamma_\Omega,\mu} \cdot \psi_\alpha \cdot 
(\phi_{\frac{\pi^2}{4\sigma^2}}- C_{L,\sigma}) \, , $$
with $\psi_\alpha=\fF(\phi_\alpha)$. It is clear from this expression that $f_{\Gamma_\Omega,\mu}$ is a
real valued Schwartz function. 
Since $|\det(C^t B)| \,  \cD^2_{\Gamma_\Omega,\mu}(x) \cdot \psi_\alpha(x) \geq 0$
for all $x\in \R^n$, while $\phi_{\frac{\pi^2}{4\sigma^2}}(x) \geq C_{L,\sigma}$ iff $\| x \|\leq \ell_{L^\vee}$,  we have
$$ f_{\Gamma_\Omega,\mu}(x) \leq 0 \ \ \ \text{ for } \ \ \  \| x \|\geq \ell_{L^\vee} $$
and $f_{\Gamma_\Omega,\mu}(x) \geq 0$ otherwise. 

\smallskip

When computing Fourier transforms, we interpret the Fourier transform $\fF(h_\Lambda)$ in the distributional sense,
so that we have
$$ \fF(h_\Lambda) = \psi_{\frac{\pi^2}{4\sigma^2}} -  C_{L,\sigma}\,  \delta_0\, , $$
with $\delta_0$ the Dirac delta distribution centered at $0$. 
The convolution product of the Dirac delta distribution $\delta_0$ with
a test function $\varphi$ leaves the test function unchanged,
$$ (\varphi\star \delta_0)(x) = \int_{\R^n} \varphi(x-u) \delta_0(u)\, du =\varphi(x)\, . $$
Thus, we obtain 
$$ \fF(f_{\Gamma_\Omega,\mu}) =|\det(C^t B)| \,  \fF(\cD^2_{\Gamma_\Omega,\mu})\star \phi_\alpha \star ( \psi_{\frac{\pi^2}{4\sigma^2}} - 
C_{L,\sigma} \delta_0) = $$
$$ |\det(C^t B)| \,  \fF(\cD^2_{\Gamma_\Omega,\mu})\star \phi_\alpha \star  \psi_{\frac{\pi^2}{4\sigma^2}} -
C_{L,\sigma} \,  |\det(C^t B)| \,  \fF(\cD^2_{\Gamma_\Omega,\mu})\star \phi_\alpha \, . $$

We have
$$ \cD^2_{\Gamma_\Omega,\mu} =( 1+ 2 \sum_{\ell\in \Gamma_\Omega} \mu_\ell e_\ell)^2 =
1+ 4 \sum_{\ell\in \Gamma_\Omega} \mu_\ell e_\ell +4 \sum_{\ell,\ell' \in \Gamma_\Omega} \mu_\ell  \mu_{\ell'} 
e_{\ell+\ell'} $$
so that the Fourier transform $\fF(\cD^2_{\Gamma_\Omega,\mu})$, also interpreted in the distributional sense, gives
$$ \fF(\cD^2_{\Gamma_\Omega,\mu}) = 1+ 4 \sum_{\ell\in \Gamma_\Omega} \mu_\ell \delta_\ell +4 \sum_{\ell,\ell' \in \Gamma_\Omega} \mu_\ell  \mu_{\ell'} 
\delta_{\ell+\ell'} \, , $$
with $\delta_{x_0}$ the Dirac delta centered at $x_0$. 

The convolution product $\fF(\cD^2_{\Gamma_\Omega,\mu})\star \phi_\alpha$ is then given by
$$ \fF(\cD^2_{\Gamma_\Omega,\mu})\star \phi_\alpha =\phi_\alpha +
4 \sum_{\ell\in \Gamma_\Omega} \mu_\ell T_\ell \phi_\alpha +4 \sum_{\ell,\ell' \in \Gamma_\Omega} \mu_\ell  \mu_{\ell'} 
T_{\ell+\ell'} \phi_\alpha\, . $$
Thus, we obtain a non-negative Fourier transform $\fF(f_{\Gamma_\Omega,\mu})(x)\geq 0$ for all $x\in \R^n$ if the
following inequality holds, for all $\ell\in \Gamma_\Omega$ and all $x\in \R^n$:
\begin{equation}\label{mainestimate}
(T_\ell \phi_\alpha \star  \psi_{\frac{\pi^2}{4\sigma^2}})(x) \geq C_{L,\sigma} \,\, T_\ell \phi_\alpha(x) \, .
\end{equation}
Since we are assuming that $L^\vee$ is a critical lattice in $\R^n$, we  have as in \eqref{CLsigmaest}
 $$ C_{L,\sigma} \leq e^{-\frac{n}{2}\, \frac{\pi}{4\sigma e}  | L^\vee |^{1/2n} }\, . $$
Then for $\alpha$ and $\sigma$ in the range \eqref{alpharange}, we obtain from Lemma~\ref{singleGaussian}
that \eqref{mainestimate} is verified.
\endproof

\medskip
\subsection{Lattice solutions}

The Cohn-Elkies functions constructed above are associated to a lattice $L\subset \R^n$ with the
property that its dual $L^\vee$ is a critical lattice, namely one whose sphere packing
density is maximal among lattices. The Voronoi algorithm provides a way to enumerate all these 
locally optimal solutions of the lattice-packing problem, by describing 
the space of lattices up to isometry in terms of positive definite quadratic forms and the identifying
the local maxima of the density function with the vertices of the Ryshkov polyhedron, \cite{Ryshkov}.
In general these local maxima will not be actual solutions of the sphere packing problem, as
the actual solution may not be achievable by a lattice. In terms of Cohn-Elkies functions, the
property that the critical lattice $L^\vee$ is also an actual solution of the sphere packing problem
is reflected in whether the additional property \eqref{CEid} can also be satisfied, as this provides a
sufficient condition for the optimality of the lattice for the sphere packing problem (see Lemma~\ref{CElatt}). 

\smallskip

In the case of the Cohn-Elkies functions of Theorem~\ref{mainthm}, this optimality condition can
be described more explicitly as follows.

\smallskip

\begin{prop}\label{CEprop4}
Setting
\begin{equation}\label{Upsilon}
 \Upsilon_{\alpha,\cF,\mu}:=1+4 \sum_{\ell\in \cF} \mu_\ell e^{-\alpha \| \ell \|^2} +4 \sum_{\ell,\ell' \in \cF} \mu_\ell  \mu_{\ell'} e^{-\alpha \| \ell+\ell' \|^2}
 \end{equation} 
we obtain that the condition 
\begin{equation}\label{prop4}
f_{\Gamma_\Omega,\mu}(0)=\frac{1}{|L^\vee|} (\fF f_{\Gamma_\Omega,\mu})(0) 
\end{equation}
is given by
\begin{equation}\label{prop4expl}
\begin{array}{c} \displaystyle{
  \Upsilon_{\frac{\alpha}{2},\Gamma_\Omega,\mu} \cdot \left( (1-e^{-\sigma \ell^2_{L^\vee_{2\sigma}}}) \left( \frac{\pi}{2\alpha}\right)^{n/2} + \frac{1}{|L^\vee|} e^{-\sigma \ell^2_{L^\vee_{2\sigma}}} \right) =} \\[3mm] \displaystyle{
\Upsilon_{\frac{4\sigma^2 \alpha}{\alpha+4\sigma^2},\Gamma_\Omega,\mu} \cdot \frac{1}{|L^\vee|} \cdot \left( \frac{4\sigma^2}{\alpha + 4 \sigma^2} \right)^{n/2} \, . } \end{array}
\end{equation}
\end{prop}

\proof
We have 
$$ f_{\Gamma_\Omega,\mu} (0) =\langle \bT_{\Gamma_\Omega,\mu} \gamma , \phi_\alpha \rangle 
\cdot h_\Lambda(0) = \langle \bT_{\Gamma_\Omega,\mu} \gamma , \phi_\alpha \rangle \,  
(1- e^{-\sigma \ell^2_{L^\vee_{2\sigma}}}) $$
$$ =(1- e^{-\sigma \ell^2_{L^\vee_{2\sigma}}}) |\det(C^t B)|\, ( \langle \gamma, \phi_\alpha \rangle +2\sum_{\ell \in \cF} \mu_\ell \langle T_\ell \gamma, \phi_\alpha \rangle ) $$
$$ =(1- e^{-\sigma \ell^2_{L^\vee_{2\sigma}}}) |\det(C^t B)|\, ( \langle\phi_\alpha,\phi_\alpha\rangle +
4 \sum_{\ell\in \Gamma_\Omega} \mu_\ell \langle T_\ell \phi_\alpha,\phi_\alpha\rangle +4 \sum_{\ell,\ell' \in \Gamma_\Omega} \mu_\ell  \mu_{\ell'} 
\langle T_{\ell+\ell'} \phi_\alpha,\phi_\alpha\rangle)\, . $$
Writing $e^{-(\alpha \| x \|^2+ \alpha \| x-\ell \|^2)} =e^{-\alpha \| \ell \|^2/2} \, e^{- 2\alpha \| x- \frac{\ell}{2} \|^2}$
we get
$$ f_{\Gamma_\Omega,\mu} (0) =(1- e^{-\sigma \ell^2_{L^\vee_{2\sigma}}}) |\det(C^t B)|\,\left( \frac{\pi}{2\alpha}\right)^{n/2}
(1+4 \sum_{\ell\in \Gamma_\Omega} \mu_\ell e^{-\frac{\alpha}{2} \| \ell \|^2} +4 \sum_{\ell,\ell' \in \Gamma_\Omega} \mu_\ell  \mu_{\ell'} e^{-\frac{\alpha}{2} \| \ell+\ell' \|^2}) \, . $$
On the other hand we have
$$ (\fF f_{\Gamma_\Omega,\mu}) (0) =|\det(C^t B)|\,  (\langle\phi_\alpha,\psi_{\frac{\pi^2}{4\sigma^2}}\rangle +
4 \sum_{\ell\in \Gamma_\Omega} \mu_\ell \langle T_\ell \phi_\alpha,\psi_{\frac{\pi^2}{4\sigma^2}}\rangle +4 \sum_{\ell,\ell' \in \Gamma_\Omega} \mu_\ell  \mu_{\ell'} 
\langle T_{\ell+\ell'} \phi_\alpha,\psi_{\frac{\pi^2}{4\sigma^2}}\rangle) $$
$$ - |\det(C^t B)|\, e^{-\sigma \ell^2_{L^\vee_{2\sigma}}}
(1+4 \sum_{\ell\in \Gamma_\Omega} \mu_\ell e^{-\frac{\alpha}{2} \| \ell \|^2} +4 \sum_{\ell,\ell' \in \Gamma_\Omega} \mu_\ell  \mu_{\ell'} e^{-\frac{\alpha}{2} \| \ell+\ell' \|^2}) \, . $$
We have
$$ \left( \frac{4\sigma^2}{\pi} \right)^{n/2} \int_{\R^n} e^{-\alpha \| x-\ell \|^2} e^{-4\sigma^2 \| x \|^2} dx =
\left( \frac{4\sigma^2}{\alpha + 4 \sigma^2} \right)^{n/2} \, e^{-\frac{4\sigma^2 \alpha}{\alpha+4\sigma^2} \| \ell \|^2} \, , $$
hence we obtain
$$ (\fF f_{\Gamma_\Omega,\mu}) (0) =|\det(C^t B)|\, \left( \frac{4\sigma^2}{\alpha + 4 \sigma^2} \right)^{n/2} \,
(1+4 \sum_{\ell\in \Gamma_\Omega} \mu_\ell e^{-\frac{4\sigma^2 \alpha}{\alpha+4\sigma^2} \| \ell \|^2} $$ $$+4 \sum_{\ell,\ell' \in \Gamma_\Omega} \mu_\ell  \mu_{\ell'} e^{-\frac{4\sigma^2 \alpha}{\alpha+4\sigma^2} \| \ell+\ell' \|^2}) $$
$$ - |\det(C^t B)|\, e^{-\sigma \ell^2_{L^\vee_{2\sigma}}}
(1+4 \sum_{\ell\in \Gamma_\Omega} \mu_\ell e^{-\frac{\alpha}{2} \| \ell \|^2} +4 \sum_{\ell,\ell' \in \Gamma_\Omega} \mu_\ell  \mu_{\ell'} e^{-\frac{\alpha}{2} \| \ell+\ell' \|^2}) \, . $$
Thus, using \eqref{Upsilon} we can write \eqref{prop4} in the form \eqref{prop4expl}, where $$|L^\vee|^{-1}=|L|=|\det(C)|.$$
\endproof

\smallskip

Note that, while we always have the vanishing at lattice points of the
function $\langle \gamma, \pi_z \phi \rangle$ by the Wexel--Raz duality, 
in general we do not have the vanishing of the $\langle \pi_\ell \gamma, \pi_z \phi \rangle$,
hence of $f_{\Gamma_\Omega,\mu}$. This vanishing occurs in the case when
this is a special Cohn-Elkies function (by Corollary~\ref{CEzeros})
that is, when the identity of Proposition~\ref{CEprop4} holds. 

\medskip

\section{From lattices to periodic sets}

The construction presented above uses essentially the fact that the sphere packing
considered is a lattice sphere packing. Since one expects that only in very few 
dimensions the optimal sphere packing will be realized by a lattice, one would like
to extend this method of construction of Cohn--Elkies functions adapted to lattices
to the case of periodic sets, which are known to approximate the maximal
density in any dimension. The Zassenhaus conjecture predicts that 
in every dimension the maximal density sphere packing can be
realized by a periodic packing.

\smallskip

\begin{defn}\label{perset}
A periodic set $\Sigma$ in $\R^n$ is a set for which there exist a finite collection
$\{ a_1, \ldots, a_N \}$ of vectors $a_i \in \R^n$ and a lattice $L \subset \R^n$
such that
\begin{equation}\label{persetSigma}
 \Sigma = \bigcup_{i=1}^N a_i + L \, . 
\end{equation} 
The lattice $L$ can be taken to be the maximal period lattice for $\Sigma$.
The {\rm size} of a periodic set is the minimal number $N$ of translations
such that the set can be represented in the form \eqref{persetSigma}.
A periodic set $\Sigma\subset \R^n$ of size $N$ is {\rm critical} if it is
a maximizer of the sphere packing density in $\R^n$ among all 
periodic packings of size at most $N$.
\end{defn}

\smallskip

The center density of a sphere packing with sphere centers places at the
points of a periodic set $\Sigma$ is given by
\begin{equation}\label{centerdenseN}
 \delta_\Sigma =\frac{N \, \ell^n_\Sigma}{2^n |L|}\, ,
\end{equation}
with the minimal length $\ell_\Sigma$, which is given by
$$ \ell_\Sigma= \min\{ \| \ell + a_i - a_j \|\,|\, \ell\in L,\,\,  i,j=1,\ldots, N\}\, . $$

\smallskip

\begin{defn}
A periodic set $\Sigma\subset \R^n$ of size $N$ is {\rm critical} if it is
a maximizer of the sphere packing density $\delta_\Sigma$ of \eqref{centerdenseN}
in $\R^n$, among all the periodic packings of size at most $N$, with fixed ratio $N/|L|$,
that is, if it maximizes $\ell_\Sigma$ among all periodic sets of size at most $N$.
\end{defn}

\smallskip

Gabor frames $\cG(\phi,\Lambda)$ where $\Lambda$ is not a lattice have been
considered in signal analysis, though a lot less is known about them than in the
lattice case. For Gabor frames with irregular and semi-regular $\Lambda\subset \R^2$
see for instance \cite{BeKuLyu}, \cite{CaChr}, \cite{GrRoSt}. 

\smallskip

Here we need to consider a special type of irregular Gabor frames, namely semiregular 
frames with $\Lambda=\Sigma \times K \subset \R^{2n}$ where $\Sigma\subset \R^n$
is a periodic set and $K\subset \R^n$ is a lattice, and with a Gaussian window function. 

\smallskip

\begin{defn}\label{Multisys}
A Gabor multisystem $\cG(\phi_1,\ldots,\phi_N, \Lambda_1, \ldots, \Lambda_N)$ is
defined as the union of the Gabor systems
$$ \cG(\phi_1,\ldots,\phi_N, \Lambda_1, \ldots, \Lambda_N):= \bigcup_{i=1}^N \cG(\phi_i, \Lambda_i)\, . $$
A multi-window Gabor system is a multisystem of the form $\cG(\phi_1,\ldots,\phi_N, \Lambda_1, \ldots, \Lambda_N)$ where $\Lambda_i=\Lambda$ for all $i=1,\ldots,N$. We write $\cG(\underline{\phi},\Lambda)$
in this case.
\end{defn}

\begin{rem}\label{multirem}{\rm 
Note that the Gabor system $\cG(\phi, \Lambda)$ with $\Lambda=\Sigma \times K$
for $\Sigma=\cup_{i=1}^N a_i + L$ is the same as the multi-window system given by the union of the
$\cG(\phi_i,\Lambda_i)$ where $\phi_i=\pi_{a_i}\phi$ and with $\Lambda_i = L\times K$
for all $i=1,\ldots, N$. }
\end{rem}

\smallskip

\begin{rem}\label{detABcond}{\rm 
Theorem~5.1 of \cite{BowRze} shows that a necessary condition for
completeness of the Gabor multi-window system $\cG(\underline{\phi}, L\times K)$, 
with $\phi_i=\pi_{a_i}\phi$ and $\phi$ a Gaussian, 
is the condition that $\det(A) \times \det(B) < N$, 
where $L= A \Z^n$ and $K= B \Z^n$ for $A,B\in \GL_n(\R)$ and
with $N$ the number of translations of the periodic set $\Sigma$. }
\end{rem}

\smallskip

In the case of a periodic set \eqref{persetSigma} we define
\begin{equation}\label{Sigmavee}
\Sigma^\vee :=\bigcup_{i=1}^N a_i + L^\vee \, ,
\end{equation}
where $L^\vee$ is the dual lattice of $L$ and the $\{ a_i \}_{i=1}^N$ are the same translations
as in $\Sigma$. We then set
\begin{equation}\label{CorrSigma}
Corr(\phi_{2\sigma}, \Sigma^\vee) := \max_{\ell,\ell'\in \Sigma^\vee, \ell\neq \ell'} \,\, | \langle \pi_{i\ell} \phi_{2\sigma}, \pi_{i\ell'} \phi_{2\sigma} \rangle | \, ,
\end{equation}
which directly generalizes \eqref{Corr2n} for lattices. We also define as before
\begin{equation}\label{CsigmaSigma}
C_{\Sigma,\sigma}:= \left(\frac{4\sigma}{\pi}\right)^{n/2}\, Corr(\phi_{2\sigma},\Sigma^\vee_{2\sigma}) \, ,
\end{equation}
with $\Sigma^\vee_{2\sigma}:=(\Sigma^\vee)_{2\sigma}=\frac{\pi}{2\sigma} \Sigma^\vee$.

\smallskip

\begin{prop}\label{mainpropSigma}
Let $\Sigma\subset \R^n$ be a periodic set of size $N$ with $\Sigma^\vee$ critical, 
and let $K\subset \R^n$ be a lattice such that the $\cG(\phi_{j,\alpha}, L\times K)$ for $j=1,\ldots, N$
are Gabor frames. Then a Cohen-Elkies functions in dimension $n$
with size $\ell_{\Sigma^\vee}$ is given by
$f_{\Gamma_\Omega,\mu}=\sum_{j=1}^N f_{j,\Gamma_\Omega,\mu}$,
with
$$ f_{j,\Gamma_\Omega,\mu} (x) :=\langle \bT_{\Gamma_\Omega,\mu} \gamma_j , \pi_{ix} \phi_{j,\alpha} \rangle \cdot h_\Sigma(x)\, , $$
$$ h_\Sigma(x)=\phi_{\frac{\pi^2}{4\sigma^2}}(x) - C_{\Lambda,\sigma}\, , $$ 
with $\alpha,\sigma$ in the range \eqref{alpharange}. 
\end{prop}

\smallskip

\proof
In general a multi-window Gabor system $\cG(\underline{\phi},L \times K)$, 
with $\Lambda\subset \R^{2n}$ a lattice, satisfies the frame condition  if
there are constants $C,C'>0$ such that
for all $f\in L^2(\R^n)$
$$ C \| f \|_{L^2(\R^n)}^2 \leq \sum_{j=1}^N \sum_{\lambda\in L\times K} |\langle f, \pi_\lambda \pi_{a_j} \phi \rangle |^2
\leq C' \| f \|_{L^2(\R^n)}^2\, . $$
Thus, if the individual Gabor systems $\cG(\phi_j, L\times K)$ satisfy the Gabor frame condition
then the multi-window system also does. 

\smallskip

Since only the translations part $\Sigma\subset \R^n$ is a periodic set, 
while the modulation part $K\subset \R^n$ is an actual lattice, we have that the
functions of the multi-window $\phi_j=\pi_{a_j}\phi= T_{a_j}\phi$ are translates
of the Gaussian $\phi$ along the translations $a_j$ of the periodic set 
$\Sigma=\cup_j a_j+L$.
 
\smallskip

Given the periodic set $\Sigma\subset \R^n$, suppose that the lattice $K\subset \R^n$ is
chosen so that the condition of Remark~\ref{detABcond} holds and the Gabor systems
$\cG(\phi_{j,\alpha}, L\times K)$, with $\phi_{j,\alpha}=T_{a_j}\phi_\alpha$ and $\phi_\alpha$ a Gaussian 
are Gabor frames for all $i=1,\ldots, N$. 

\smallskip 

Then proceeding as in Proposition~\ref{approxdual}, with $\Omega >> \| a_j \|$ for
all $i=1,\ldots, N$, we obtain approximate Wexel-Raz duals $\gamma_j$ for each
$\cG(\phi_{j,\alpha}, L\times K)$. Note that the sets $\Gamma_\Omega$ and the
multiplicity function $\mu: \Gamma_\Omega \to \N$ are unchanged, and we have
approximate Wexel--Raz duals for the Gabor systems $\cG(\phi_{j,\alpha}, L\times K)$
of the form
$$ \gamma_j =|\det(A^t B)| \, (\phi_{j,\alpha}+ 2\sum_{\ell\in \Gamma_\Omega} \mu_\ell T_\ell \phi_{j,\alpha}) =
|\det(A^t B)| \,  \bT_{\Gamma_\Omega,\mu} \phi_{j,\alpha} \, . $$

\smallskip

As in \eqref{corr2sigma2} we have
\begin{equation}\label{corr2sigma2Sigma}
  |\langle \pi_{i\ell} \phi_{2\sigma}, \pi_{i\ell'} \phi_{2\sigma} \rangle|=\left( \frac{\pi}{4\sigma}\right)^{n/2} \, e^{-\frac{\pi^2}{4\sigma} \| \ell-\ell' \|^2} \, , 
\end{equation}  
hence we obtain
\begin{equation}\label{ellSigmaCorrSigma}
C_{\Sigma,\sigma} = \max_{\ell,\ell'\in \Sigma^\vee, \ell\neq \ell'} 
e^{-\frac{\pi^2}{4\sigma} \| \ell-\ell' \|^2} = e^{-\frac{\pi^2}{4\sigma} \ell^2_{\Sigma^\vee}} \, .
\end{equation}

\smallskip

We can assume without loss of generality that the lattices $L\times K$
we consider all have $L\subset \R^n$ with $|L|=1$, with conditions 
such as Remark~\ref{detABcond} for the Gabor frames property formulated 
as conditions on the choice of the auxiliary lattice $K\subset \R^n$.

\smallskip

Assuming that the periodic set $\Sigma^\vee$ is critical, we have 
$$ \ell^2_{\Sigma^\vee} \geq \gamma_n \geq \frac{n}{2\pi e}\, , $$
since $\Sigma^\vee$ maximizes density among all packings
by periodic sets of size at most $N$, hence its density is not
worse than the optimal density among lattices. We use
the same lower bound for the Hermite constant $\gamma_n$ as in
Corollary~\ref{LCn}. Thus, we obtain
\begin{equation}\label{boundCSigma}
C_{\Sigma,\sigma} \leq e^{-n \frac{\pi}{8\sigma e}},
\end{equation}
as in \eqref{CLsigmaest}. We can then use the
same estimates of Lemma~\ref{singleGaussian} and Remark~\ref{XiL}
(with $|L|=|L^\vee|=1$). 

\smallskip

Lemma~\ref{singleGaussian} holds for each $g_{j,\alpha}=T_\ell \phi_{j,\alpha}$, hence
the same argument used in Theorem~\ref{mainthm} shows that the functions
$$ f_{j,\Gamma_\Omega,\mu} (x) :=\langle \bT_{\Gamma_\Omega,\mu} \gamma_j , \pi_{ix} \phi_{j,\alpha} \rangle \cdot h_\Sigma(x)\, , $$
with 
$$ h_\Sigma(x)=\phi_{\frac{\pi^2}{4\sigma^2}}(x) - C_{\Lambda,\sigma}\, , $$ 
with the parameters $\alpha,\sigma$ in the same range as in Theorem~\ref{mainthm},
are Cohen-Elkies functions in dimension $n$
with size $\ell_{\Sigma^\vee}$, as in Definition~\ref{CEfunction}, hence so is their sum 
$f_{\Gamma_\Omega,\mu}=\sum_{j=1}^N f_{j,\Gamma_\Omega,\mu}$. 
\endproof

\bigskip

\subsection*{Acknowledgment} The second author is supported by NSF grant DMS-2104330.


\begin{thebibliography}{99}

\bibitem{BaCaHeLa} R.~Balan, P.G.~Casazza, C.~Heil, Z.~Landau, {\em Density, 
overcompleteness, and localization of frames II: Gabor frames}, J. Fourier Anal.
Appl., 12 (2006) N.3, 307--344.

\bibitem{Bale} R.A.~Bale, J.P.~Grossman, G.F.~Margrave, M.P.~Lamoureux,
{\em Multidimensional partitions of unity and Gaussian terrains}, CREWES Research Report, Volume 14 (2002).

\bibitem{BeKuLyu} Yu.~Belov, A.~Kulikov, Yu.~Lyubarskii, {\it Irregular Gabor frames of Cauchy kernels},
arXiv:2104.01121.

\bibitem{BowRze} M.~Bownik, Z.~Rzeszotnik, {\it The spectral function of shift-invariant spaces}, 
Michigan Math. J. 51 (2003), no. 2, 387--414.

\bibitem{CaChr} P.G.~Casazza, O.~Christensen, {\it Gabor frames over irregular lattices}, 
Advances in Computational Mathematics 18 (2003) 329--344.

\bibitem{ChrKim} O.~Christensen, R.Y.~Kim, {\em Pairs of explicitly given dual Gabor frames in $L^2(\R^d)$},
The Journal of Fourier Analysis and Applications, 
Vol.~12 (2006), N.~3, 243--255. 

\bibitem{CoEl} H.~Cohn, N.~Elkies, {\it New upper bounds on sphere packings, I},
Ann. of Math. 157 (2003) 689--714.

\bibitem{CoKu} H.~Cohn, A.~Kumar, {\it
Universally optimal distribution of points on spheres}, J. Amer. Math. Soc. 20 (2007) 99--148.

\bibitem{CKMRV} H.~Cohn, A.~Kumar, S.D.~Miller, D.~Radchenko, M.~Viazovska,
{\it The sphere packing problem in dimension 24}, 
Ann. of Math. (2) 185 (2017), no. 3, 1017--1033.

\bibitem{Groch} K.~Gr\"ochenig, {\it  Foundations of Time--Frequency Analysis},  Birkh\"auser,  2001.

\bibitem{GroKop} K.~Gr\"ochenig, S.~Koppensteiner, {\it Gabor Frames: Characterizations and Coarse Structure},
in ``New Trends in Applied Harmonic Analysis, Volume 2", pp.~93--120, Birkh\"auser, 2019. 

\bibitem{GrRoSt} K.~Gr\"ochenig, J.L.~Romero, J.~St\"ockler,
{\it Sampling Theorems for Shift-invariant Spaces, Gabor Frames, and Totally Positive Functions},
Inventiones Mathematicae, 211(2018) N.3, 1119--1148.

\bibitem{Jan2} A.J.E.M.~Janssen, {\em Signal analytic proofs of two basic results on lattice expansions}, 
Appl. Comput. Harmon. Anal., Vol.1 (1994) N.4, 350--354.

\bibitem{Papush} D.E.~Papush, {\it Interpolation with discrete sets in $\C^\ell$}, J. Soviet
Math., Vol.59 (1992) N.1, 666--674.

\bibitem{Ryshkov} S.S.~Ryshkov, {\em The polyhedron $\mu(m)$ and certain extremal problems of the geometry of numbers}, 
Soviet Math. Dokl, Vol. 11 (1970), 1240--1244.

\bibitem{StrohHeath} T.~Strohmer, R.W.~Heath Jr., {\em Grassmannian frames with 
applications to coding and communication}, Appl. Comput. Harmon. Anal. 14 (2003) 257--275. 

\bibitem{Viaz} M.~Viazovska, {\it The sphere packing problem in dimension 8},  
Ann. of Math. (2) 185 (2017), no. 3, 991--1015. 


\end{thebibliography}
\end{document}